\documentclass[a4paper,12pt]{article}
\usepackage{amssymb}

\newcommand{\R}{\mathbb{R}}

\newcommand{\N}{\mathbb{N}}

\newcommand{\beq}{\begin{equation} }
\newcommand{\eqq}{\end{equation} }
\newcommand{\cuad}{{\sqcap\kern-.68em\sqcup}}

\newcommand{\equ}[1]{(\ref{#1})}

\newtheorem{teo}{Theorem}[section]

\newtheorem{proposition}{Proposition}[section]

\newtheorem{lemma}{Lemma}[section]

\newtheorem{remark}{Remark}[section]
\newcommand{\bremark}{\begin{remark} \em}
\newcommand{\eremark}{\end{remark} }

\def\beeq{\begin{equation}}
\def\eeq{\end{equation}}
\newcommand{\begeqaet}{\begin{eqnarray*}}
\newcommand{\eneqaet}{\end{eqnarray*}}

\begin{document}
\begin{center}{\bf \Large   Qualitative properties of positive solutions \medskip

for mixed integro-differential equations }\medskip

\bigskip

{Patricio Felmer\ \  and\ \  Ying Wang
\\~\\}

Departamento de Ingenier\'{\i}a  Matem\'atica and Centro de
Modelamiento Matem\'atico,
UMR2071 CNRS-UChile,
 Universidad de Chile\\
{\sl   ( pfelmer@dim.uchile.cl   and   yingwang00@126.com ) }

\begin{abstract}
This paper is concerned with the  qualitative properties  of  the solutions of
mixed integro-differential equation
\begin{equation}\label{eq 1}
\left\{ \arraycolsep=1pt
\begin{array}{lll}
 (-\Delta)_x^{\alpha} u+(-\Delta)_y u+u=f(u)\quad \ \  {\rm in}\ \  \R^N\times\R^M,\\[2mm]
 u>0\ \ {\rm{in}}\  \R^N\times\R^M,\ \ \quad
 \lim_{|(x,y)|\to+\infty}u(x,y)=0,
\end{array}
\right.
\end{equation}
with $N\ge 1$, $M\ge 1$ and $\alpha\in (0,1)$. We study decay and symmetry properties of the solutions to this equation. Difficulties arise due to the   mixed character of the integro-differential operators. Here, a crucial role is played by a version of the Hopf's Lemma we prove in our setting.  In studying the decay, we
 construct appropriate super and sub solutions and we use the moving planes method to prove the symmetry properties.
\end{abstract}
\end{center}

  {\small \noindent {\bf Key words}:  Integro-differential equation, Hopf's Lemma,  Decay,  Symmetry.\\
 \noindent {\bf MSC2010}:  35R11, 35B06, 35B40, 35B50.}

\setcounter{equation}{0}
\section{Introduction}

The study of qualitative properties of positive solutions to semi-linear elliptic equations in $\R^N$ has been
the concern of numerous authors along the last several decades. The asymptotic behavior of the solution at infinity, the actual rate of decay and symmetry properties have been the most studied qualitative properties for these equations. It was the seminal work by Gidas, Ni and Nirenberg \cite{GNN2} that settled these two main qualitative properties for the semi-linear elliptic equation
\begin{equation}\label{eq 2.1wy}
\left\{ \arraycolsep=1pt
\begin{array}{lll}
 -\Delta u+u=f(u)\quad {\rm in}\quad \R^M,\\[2mm]
 u>0\ \ {\rm{in}}\ \  \R^M,\ \ \quad
 \lim_{|y|\to+\infty}u(y)=0,
\end{array}
\right.
\end{equation}
 when the non-linearity is merely Lipschitz continuous, super-linear at the zero, in the sense that
  \begin{equation}\label{1.2}
f(s)=O(s^p) \quad {\rm as}\ s\to 0,
\end{equation}
 for some $p>1$, and  $M\ge 3$. Gidas, Ni and Nirenberg proved that  the solutions of \equ{eq 2.1wy} are  radially symmetric and they satisfy the precise decay estimate
\begin{equation}\label{1.1jx1}
 \lim_{|y|\to+\infty}u(y)e^{|y|}|y|^{\frac{M-1}{2}}=c,
\end{equation}
for certain constant $c>0$. After this work, many authors extended the results in various directions, generalizing the non-linearity, the elliptic operator or the hypotheses on the solutions. Out of the very many contributions in this direction we mention here only a few:  Berestycki and Lions \cite{BL}, Berestycki and Nirenberg \cite{BN}, Brock \cite{Br}, Busca and Felmer \cite{BF}, Cort\'azar, Elgueta and Felmer \cite{CEF2}, Da Lio and Sirakov \cite{DS}, Dolbeault and Felmer \cite{DF}, Gui \cite{Gui}, Kwong \cite{KK}, Li and Ni \cite{LN} and  Pacella and Ramaswamy \cite{PR}.

Recently, much attention has been given to the study of elliptic equations of fractional order. In this direction, Felmer, Quaas and Tan in \cite{FQT} studied the problem
\begin{equation}\label{eq 2.1wy000}
\left\{ \arraycolsep=1pt
\begin{array}{lll}
 (-\Delta)^\alpha u+u=f(u)\quad {\rm in}\quad \R^N,\\[2mm]
 u>0\ \ {\rm{in}}\ \  \R^N,\ \ \quad
 \lim_{|x|\to+\infty}u(x)=0.
\end{array}
\right.
\end{equation}
They proved existence and regularity of positive solutions, and also decay and symmetry results. Precisely, it was proved that the solutions $u$ of  (\ref{eq 2.1wy000}) satisfy
\begin{equation}\label{1.1}
\frac{c^{-1}}{|x|^{N+2\alpha}}\le u(x)\le \frac{c}{|x|^{N+2\alpha}},\qquad |x|\ge1,
\end{equation}
for some $c>1$, when $f$ is superlinear at $0$ in the sense that
$$
\lim_{s\to 0}\frac{f(s)}{s}=0.
$$
The radial symmetry of the solutions of   (\ref{eq 2.1wy000})  is derived
by using the moving planes method in integral form developed in \cite{CLO2,L}, assuming further that
$f\in C^1(\R),$ it is increasing and there exists $\tau>0$ such that
 \begin{equation}\label{eq fqt1}
 \lim_{s\to0}\frac{f'(s)}{s^{\tau}}=0.
 \end{equation}
This symmetry result was generalized by the authors in  \cite{FW}, using an appropriate truncation argument together with the moving planes method with ideas developed in \cite{LCM}.
We refer to some other papers with more discussions on qualitative properties of solutions to
fractional elliptic problems as Cabr\'e and Sire \cite{CS5}, Caffarelli and Silvestre \cite{CS}, Chen, Li and Ou \cite{CLO2}, Barles, Chasseigne, Ciomaga and Imbert \cite{BCCI3}, Dipierro Palatucci, Valdinoci \cite{DPV}, Li \cite{L}, Quaas and Xia \cite{QX},
Ros-Oton and Serra \cite{RS} and Sire and Valdinoci \cite{S1}.

Both operators, the laplacian and the fractional laplacian, are particular cases of a general class of elliptic operators connected to backward stochastic differential equations associated to Brownian and Levy-It\^{o} processes, see for example Barles,  Buckdahn and Pardoux \cite{BBP}, Benth, Karlsen and Reikvam \cite{BKR} and Pham \cite{P}. Recently,       Barles, Chasseigne, Ciomaga and Imbert in   \cite{BCCI,BCCI2} and Ciomaga in \cite{BCCI3} considered
 the  existence and  regularity of solutions for  equations involving  mixed  integro-differential operators belonging  to the  general class of backward stochastic differential equations mentioned above. A particular case of   elliptic integro-differential operator  of mixed type is the one considering the laplacian in some of the variables  and the fractional laplacian in the others, modeling diffusion sensible to the direction. In view of \equ{eq 2.1wy} and \equ{eq 2.1wy000} we may write  similarly
\begin{equation}\label{eq wyq1}
\left\{ \arraycolsep=1pt
\begin{array}{lll}
 (-\Delta)_x^{\alpha} u+(-\Delta)_y u+u=f(u),\ \ \ \  (x,y)\in\R^N\times\R^M,\\[2mm]
 u>0\ \ {\rm{in}}\  \R^N\times\R^M,\ \ \quad
 \lim_{|(x,y)|\to+\infty}u(x,y)=0,
\end{array}
\right.
\end{equation}
where $N\geq1$, $M\geq1$. The operator  $(-\Delta)_y $ denotes the usual laplacian with respect to $y$, while
$(-\Delta)_x^\alpha $ denotes the fractional laplacian of exponent $\alpha\in(0,1)$ with respect to  $x$, i.e.
  \begin{equation}\label{DELTA}
  (-\Delta)_x^\alpha u(x,y)=\int_{\R^N}\frac{u(x,y)-u(z,y)}{|x-z|^{N+2\alpha}}dz,
  \end{equation}
for all $(x,y)\in\R^N\times\R^M$. Here the integral is understood in the principal value sense.

In view of the known results on decay and symmetry for solutions of equations \equ{eq 2.1wy} and \equ{eq 2.1wy000} just described above, it is interesting to ask if these results still hold for solutions of the equation of  mixed type \equ{eq wyq1}, where the elliptic operator represents diffusion depending on the direction in space. Regarding the asymptotic decay of solution at infinity, the question is interesting since a proper mix of the two variables should be obtained for the decay estimates. The natural way to estimate the decay is through the construction of super and sub solutions  involving the fundamental solution of the elliptic operator, which in this case is singular in $\R^N\times \{0\}$. Moreover, the solution of  \equ{eq wyq1} cannot be radially symmetric, so this property cannot be used to estimate the decay. On the other hand, regarding radial  symmetry, we may still have symmetry in $x$ and $y$, but the moving planes method would require an adequate version of the Hopf's Lemma, that we prove here.

Our  first theorem concerns  the decay of  solutions for (\ref{eq wyq1}) with general nonlinearity and it states  as follows.
\begin{teo}\label{teo r1}
Let   $\alpha\in(0,1)$, $N,M\in \N$, $N\geq1$ and $M\geq1$ and let us assume that the function $f:(0,+\infty)\to\R$ is continuous and it satisfies
\begin{equation}\label{eq 011j}
  -\infty<B:=\liminf_{v\to 0^+}\frac{f(v)}v\leq A:=\limsup_{v\to 0^+}\frac{f(v)}v<1.
\end{equation}
Let $u$  be a positive classical solution of  problem
(\ref{eq wyq1}),  then for any $\epsilon>0$ small,  there exists $C_\epsilon>1$  such that
for any $(x,y)\in\R^N\times\R^M$,
 \begin{equation}\label{eq1.2}
 C_\epsilon^{-1}(1+|x|)^{-N-2\alpha} e^{-\theta_2|y|}
  \le u(x,y)\le C_\epsilon(1+|x|)^{-N-2\alpha} e^{-\theta_1|y|},
 \end{equation}
 where
 \begin{equation}\label{theta12}
 \theta_1=\sqrt{1-A}-\epsilon\quad\mbox{ and}\quad \theta_2=\sqrt{1-B}+\epsilon.
 \end{equation}
\end{teo}

When we compare estimate \equ{eq1.2}  with  \equ{1.1jx1} for $N=0$, we first observe that in ours  an exponential decay is obtained, but with a constant $C_\epsilon$ depending on $\epsilon$, which is a parameter  controlling the rate of exponential decay. This is more clear when $A=B=0$. On the other hand we are making much more general assumptions on $f$ and, in particular, we are not making any assumption on the radial symmetry of the solution, which is crucial in proving  \equ{1.1jx1}. We do not know of a decay estimate better than
 \begin{equation}\label{eq1.2n}
 C_\epsilon^{-1} e^{-\theta_2|y|}
  \le u(y)\le C_\epsilon e^{-\theta_1|y|},\qquad y\in \R^M,
 \end{equation}
for solutions of \equ{eq 2.1wy} under assumption \equ{eq 011j} for $f$,  and where radial symmetry of the solutions is not available, like in a case where $f$ may depend on $y$.
On the other hand, when $M=0$, we recover  \equ{1.1} from \equ{eq1.2}.
For the  proof of  the decay estimate (\ref{eq1.2}) we construct suitable super and sub solutions and we use comparison principle with a version of Hopf's lemma.

When we assume further hypothesis  we can get sharper estimates for the decay of the solutions of equation (\ref{eq wyq1}).
Precisely, we have the following result:
\begin{teo}\label{teo jxcwr1}
Assume that  $\alpha\in(0,1)$, $N\geq1$, $M\geq5$ and the non-linearity  $f:(0,+\infty)\to \R$ is non-negative and it  satisfies \equ{1.2}.
 Let $u$  be a positive classical solution of  (\ref{eq wyq1}),
then  there exists a constant $c>1$  such that
for all $(x,y)\in\R^N\times\R^M$,
\begin{equation}\label{jx220}
 \frac{1}{c}{\rho}(x,y)\leq u(x,y)\leq  c{\rho}(x,y)(1+|y|)^{\frac{1}{2}},
\end{equation}
where  the function $\rho$ is defined as
\begin{equation}\label{ro}
{\rho}(x,y)=\min\{ \frac{1}{(1+|x|)^{N+2\alpha}}\, ,\, e^{-|y|}|y|^{-\frac{N}{2\alpha}-\frac M2}\,,\,  \frac{ e^{-|y|}|y|^{1-\frac{M}{2}}}{(1+|x|)^{N+2\alpha}} \}.
\end{equation}
\end{teo}

We notice that this theorem gives the expected exponential decay  for positive solutions, as suggested by (\ref{1.1jx1}), assuming  the dimension of the space satisfies  $M\ge 5$. Moreover, it gives the expected polynomial correction for the lower bound with a gap in the power for the upper bound.  This theorem is proved under the assumption  \equ{1.2} on the non-linearity, constructing super and sub solutions devised upon the  fundamental solution of $(-\Delta)_x^{\alpha} +(-\Delta)_y +id$. In our argument, a crucial role is played by the estimate already obtained in Theorem 1.1. Since the fundamental solution of   $(-\Delta)_x^{\alpha} +(-\Delta)_y +id$ has $\R^N\times\{0\}$ as  singular  set, we cannot use  the method in \cite{GNN2} in order to derive our estimate. Moreover, some other arguments in \cite{GNN2} cannot be used either because the solutions of   (\ref{eq wyq1}) are     not radial,
 since the differential operator is not radially invariant and there are no solutions depending only on one of the $x$ or $y$ variables, as can be seen from   \equ{jx220},

Even though solutions of   (\ref{eq wyq1})  are not radially symmetric, we can prove partial symmetry in each of the variables $x$ and $y$ and this is the content of our third theorem.
\begin{teo}\label{teo r1jt}
Assume that $\alpha\in(0,1)$, $N\geq1$, $M\geq1$ and the function  $f:(0,+\infty)\to\R$ is locally Lipschitz and it satisfies (\ref{eq 011j}). Moreover, we assume that  $f$ also satisfies
\begin{itemize}
\item[$(F)\ $]
 there exist $u_0>0$, $\gamma>\frac{N}{N+M}\cdot\frac{2\alpha}{N+2\alpha}$ and $\bar c>0$ such that
  \begin{equation}\label{eq rq2}
 \frac{f(v)-f(u)}{v-u}\leq \bar c v^\gamma\ \ \ \ {\rm for\ all } \quad 0<u<v<u_0.
 \end{equation}
\end{itemize}
Then, every   positive classical solution $u$ of equation
(\ref{eq wyq1}) satisfies
$$u(x,y)=u(r,s)$$
and $u(r,s)$ is strictly decreasing in $r$ and $s$,
where $ r=|x|$ and $ s=|y|$.
\end{teo}

When $N=0$, we see that  assumption $(F)$ implies  $\gamma>0$ and  (\ref{eq rq2})  coincides with the assumption considered  in \cite{LCM}. When $M=0$,  assumption $(F)$ implies that  $\gamma>\frac{2\alpha}{N+2\alpha}$ and it  coincides with the assumption considered in \cite{FW},  when the solutions is assumed to decay as a power $N+2\alpha$  at infinity.
We remark that  the operator $ (-\Delta)_x^{\alpha} +(-\Delta)_y$  is a combination of two operators with  different differential orders in $x-$variable and $y-$variable, and this produced a combined polynomial-exponential decay and does not allow for radial symmetry, but only  partial symmetry as stated in Theorem 1.3.

The proof of Theorem 1.3 is based on the moving planes method as developed in \cite{FW,LCM}. In these arguments, the  strong maximum principle plays a crucial role and it is  available for the laplacian and for the fractional laplacian. However, in the case of our mixed integro-differential operator some difficulties arise and we overcome them with  a version of the   Hopf's Lemma.


The rest of the paper is organized as follows. In Section \S 2, we introduce a version of the  Hopf's Lemma and a strong maximum principle.
In Section \S3, we  prove the decay of solutions
as in Theorem \ref{teo r1} and Theorem \ref{teo jxcwr1} by  constructing suitable super and sub solutions.
Section \S4  is devoted to prove symmetry results presented in Theorem 1.3.

 \setcounter{equation}{0}
\section{Preliminaries}
This section is devoted to study the Strong
Maximum Principle for mixed integro-differential operators as in equation \equ{eq wyq1}. To this end, we prove first a suitable form of the Hopf's Lemma.

However, before to go to this, we recall  some basic properties of the Sobolev embeddings.
If we denote the Sobolev spaces
$$H(\R^{N+M})=\{w\in L^2(\R^{N+M}) |   \int_{\R^M}\int_{\R^N}(|\xi_1|^{2\alpha}+|\xi_2|^2+1)|\hat{w}(\xi_1,\xi_2)|^2d\xi_1d\xi_2<\infty  \}$$
and
$$H^\alpha(\R^{N+M})=\{w\in L^2(\R^{N+M})\ |\   \int_{\R^{N+M}}(|\xi|^{2\alpha}+1)|\hat{w}(\xi)|^2d\xi<\infty  \},$$
with norms
$$\|w\|_H=( \int_{\R^M}\int_{\R^N}(|\xi_1|^{2\alpha}+|\xi_2|^2+1)|\hat{w}(\xi_1,\xi_2)|^2d\xi_1d\xi_2)^\frac{1}{2}$$
and
$$\|w\|_{H^\alpha}=( \int_{\R^{N+M}}(|\xi|^{2\alpha}+1)|\hat{w}(\xi)|^2d\xi)^\frac{1}{2},$$
respectively, then it is not difficult to see that the following proposition holds.
\begin{proposition} \label{lemma 1}
For $\alpha\in(0,1)$, we have that
$$H(\R^{N+M})\subset H^\alpha(\R^{N+M})\subset L^p(\R^{N+M}),$$
where the first inclusion is continuous and the second inclusion is continuous if $1\leq p\leq \frac{2(N+M)}{N+M-2\alpha}$. Moreover,
$$H(\R^{N+M})\subset  L^p_{loc}(\R^{N+M})$$
 is compact if $1\leq
p<\frac{2(N+M)}{N+M-2\alpha}$.
\end{proposition}

We devote the rest of this section to prove the  Strong Maximum Principle in our context and to this end,  we start with versions of the
 Maximum Principle and the Hopf's Lemma. In what follows, given
 $\Omega$ an open subset in $\R^N\times\R^M$, we define its closed cylindrical extension in the direction $x$ as
\begin{eqnarray*}
\tilde\Omega=\{(x,y) \in \R^N\times\R^M:\  \exists\  x'\in\R^N\
{\rm{s.t.}}\ (x',y)\in\bar\Omega\}.
\end{eqnarray*}
 Given a function  $h$ defined in an appropriate domain, we consider the mixed integro-differential operator
$$\mathcal{L}w(x,y)=(-\Delta)^\alpha_xw(x,y)+(-\Delta)_yw(x,y)+h(x,y)w(x,y).$$

\begin{lemma}\label{teo MP}
Assume that $\Omega$ is an open domain of $\R^N\times\R^M$ and the function $h:\Omega\to\R$ satisfies  $h\ge0$ in $\Omega$. If the function $w\in C(\bar \Omega)\cap
L^\infty(\tilde\Omega)$ satisfies
\begin{equation}\label{eq w1}
\left\{ \arraycolsep=1pt
\begin{array}{ll}
 \mathcal{L}w\ge 0\quad {\rm in}\ \ \Omega,\quad
 w\geq0\quad  {\rm in}\ \ \tilde\Omega\setminus\Omega,\\[2mm]
\liminf_{(x,y)\in\Omega,|(x,y)|\to\infty}w(x,y)\geq0
\end{array}
\right.
\end{equation}
then
$w\ge0$ in $\tilde{\Omega}.$
\end{lemma}
{\bf Proof.} If not, we may assume that there exists some $(x_0,y_0)\in\Omega$ such
that
$$w(x_0,y_0)=\min_{(x,y)\in\tilde\Omega} w(x,y)<0.$$
Then
\begin{eqnarray*}
(-\Delta)_x^\alpha w(x_0,y_0)=\int_{
\R^N}\frac{w(x_0,y_0)-w(z,y_0)}{|x_0-z|^{N+2\alpha}}dz
<0
\end{eqnarray*}
and
$$(-\Delta)_y w(x_0,y_0)\le0$$
and then, since $h$ is non-negative we have $\mathcal{L}w(x_0,y_0)<0,$
which contradicts  (\ref{eq w1}),
 completing the proof. \hfill$\Box$\\

It what follows we
prove a version of the Hopf's Lemma and for this purpose we need to give some conditions to the boundary of the domain where the function is defined. We say that the domain
$\Omega\subset\R^N\times\R^M$ satisfies \emph{interior cylinder
condition at $(x_0,y_0)\in\partial\Omega$} if there exist $r>0$ and
$\tilde y\in \R^M$ such that $O_r=B_r^N(x_0)\times B_r^M(\tilde y)$ satisfies
\begin{equation}\label{eq ms2}
O_r\subset\Omega\quad
{\rm{and}}\quad (x_0,y_0)\in\partial O_r,
\end{equation}
 where  $B_r^N(x_0)=\{x\in
\R^N: |x-x_0|<r\}$ and $B_r^M(\tilde y)=\{y\in \R^M: |y-\tilde
y|<r\}$ and,  obviously $|\tilde y-y_0|=r$. We define also
\begin{equation}\label{eq ms1}
D=\{(x,y)\in O_r: |x-x_0|<\frac{r}{2},\ |y-\tilde y|>\frac{r}{2}\}.
\end{equation}

\begin{lemma}\label{hopf}[Hopf's Lemma]\
Let $\Omega$ be an open set satisfying interior cylinder condition
at $(x_0,y_0)\in\partial\Omega$. Assume that $h\in L^\infty(D)$
and $w\in C(\bar \Omega)\cap\L^\infty(\tilde \Omega)$  satisfies
$$\mathcal{L}w\ge 0\quad  {\rm in}\ \Omega $$
and $$0=w(x_0,y_0)<w(x,y),\quad \forall(x,y)\in \Omega.$$
Further assume that for  $r>0$ be given in (\ref{eq ms1}) and for any $(x,y)\in D$ we have
\begin{equation}\label{eq w2}
\int_{\R^N\setminus
B_r^N(x_0)}\frac{w(z,y)}{|x-z|^{N+2\alpha}}dz\ge0.
\end{equation}
Then
\begin{equation}\label{eq ms3}
\limsup_{s\to 0^+}\frac{w(x_0,y_0)-w(x_0,y_0+s\tilde y)}s<0,
\end{equation}
moreover, if the limit exists, then
\begin{equation}\label{eq w8}
 \frac{\partial w}{\partial n}(x_0,y_0)<0,
\end{equation}
where  $ n$ is the unit  exterior normal vector of $\Omega$ at the
point $(x_0,y_0)$.
\end{lemma}
{\bf Proof.}
 Let us define
\begin{equation}\label{eq w5}
\varphi_M(y)=e^{-\beta |y-\tilde y|^2}-e^{-\beta r^2},\quad y\in \bar
B_r^M(\tilde y),
\end{equation}
where $\beta>0$ will be chosen later.
By direct computation, we have that
\begin{equation}\label{eq w3}
-\Delta\varphi_M(y)=(2M\beta-4\beta^2|y-\tilde y|^2)e^{-\beta|y-\tilde y|^2}.
\end{equation}
Next we consider the function
$$ v(x,y)=\varphi_N(x)\varphi_M(y),\quad (x,y)\in \tilde O_r,$$
where $\varphi_N$ is the first eigenfunction of Dirichlet problem
\begin{equation}\label{eq w4}
\left\{ \arraycolsep=1pt
\begin{array}{lll}
(-\Delta)^\alpha \varphi_N(x)=\lambda_1 \varphi_N(x) ,\ \ \ \ &
x\in B_{r/2}^N(x_0),\\[2mm]
\varphi_N(x)=0,& x\in\R^N\setminus B_{r/2}^N(x_0),
\end{array}
\right.
\end{equation}
where $\varphi_N$ is positive and bounded in $B_{r/2}^N(x_0)$ and  the first eigenvalue $\lambda_1,$ is positive, see Propositions 9 and 4 in \cite{SV} and \cite{SV1}, respectively.

For $(x,y)\in D$, by (\ref{eq w3}) and (\ref{eq w4}), we obtain that
\begin{eqnarray*}
\mathcal{L}v(x,y)&=&
\varphi_M(y)(-\Delta)^\alpha \varphi_N(x)+\varphi_N(x)(-\Delta\varphi_M(y))+h(x,y)\varphi_N(x)\varphi_M(y)
\\&=&\varphi_N(x)[\lambda_1\varphi_M(y)+(2M\beta-4\beta^2|y-\tilde y|^2)e^{-\beta|y-\tilde y|^2}+h(x,y)\varphi_M(y)]
\\&\le&\varphi_N(x)e^{-\beta|y-\tilde y|^2}(\lambda_1+2M\beta-\beta^2r^2+\|h\|_{L^\infty( D)}),
\end{eqnarray*}
where the last inequality holds by the fact that $0\leq\varphi_M(y)<e^{-\beta |y-\tilde y|^2}$ and $|y-\tilde y|>r/2$ in $D$.
Let us choose $\beta>0$ big enough such that
\begin{equation}\label{eq w6}
\mathcal{L}v \le0\quad {\rm in} \ D.
\end{equation}
On the other hand, since $\varphi_N(x)=0$ for $|x-x_0|\geq r/2$ and
$\varphi_M(y)=0$ for $|y-\tilde y|=r$, it is obvious
 that  $v=0$ in $A_1 \cup A_2$ where $A_1=\{(x,y)\in\tilde D:  |x-x_0|\ge r/2\}$
and $A_2=\{(x,y)\in\bar D: |y-\tilde y|=r\}$. If we define the set
$A_3:=\{(x,y)\in\bar D: |y-\tilde y|=r/2\}$, we see that $\tilde D\setminus D=A_1\cup A_2\cup A_3$. We also observe that $v$ is a bounded function in $\tilde O_r$.

Next we prove  (\ref{eq ms3}) assuming $h\ge0$. Defining
\begin{equation}\label{eq w10}
W(x,y)=
\left\{ \arraycolsep=1pt
\begin{array}{lll}
w(x,y),\quad &
(x,y)\in \bar O_r,\\[2mm]
0,\quad& (x,y)\in\tilde O_r\setminus \bar O_r
\end{array}
\right.
\end{equation}
and using  (\ref{eq w2}),  we have that for any $(x,y)\in D$,
\begin{eqnarray*}
\mathcal{L}W(x,y) =\mathcal{L}w(x,y)+\int_{\R^N\setminus
B_r^N(x_0)}\frac{w(z,y)}{|x-z|^{N+2\alpha}}dz \ge0.
\end{eqnarray*}
Combining with (\ref{eq w6}), we have that, for every $\epsilon>0$
\begin{equation}\label{eq w11}
\mathcal{L}(W-\epsilon v)\ge0 \quad {\rm{in}}\ D.
\end{equation}
Since $v$ is bounded in $\tilde O_r$, the set $A_3$ is a compact subset of $O_r$ and
$w>0$ in $O_r$, then there exists $\epsilon>0$ small such that
$$W=w\ge\epsilon v\quad {\rm{in}}\ A_3.$$
Since $v=0$ in $A_1\cup A_2$, $w\geq0$ in $\bar O_r$ and (\ref{eq w10}), we have  $W\geq0=\epsilon v$ in $A_1\cup A_2$.
Consequently,
$$W-\epsilon v\geq0 \quad {\rm{in}}\ \ \tilde D\setminus D.$$
Then we can use Lemma \ref{teo MP}, recalling that  $h\geq0$ to obtain that
$$W-\epsilon v\ge0 \quad {\rm{in}}\  D.$$
In view of the definition of $W$, since $D\subset \bar O_r$, we find that  $w-\epsilon
v\ge0$ in $D$ and noticing that  $w(x_0,y_0)=v(x_0,y_0)=0$ we obtain that
$$\frac{w(x_0,y_0)-w(x_0,y_0+s\tilde y)}s\le \epsilon\frac{v(x_0,y_0)-v(x_0,y_0+s\tilde y)}s,$$
 for all $s\in(0,r/2)$. Thus, we have
\begin{eqnarray*}
\limsup_{s\to 0^+}\frac{w(x_0,y_0)-w(x_0,y_0+s\tilde y)}s &\le&
\epsilon\lim_{s\to 0^+}\frac{v(x_0,y_0)-v(x_0,y_0+s\tilde y)}s
\\&=&\epsilon \varphi_N(x_0) \lim_{s\to 0^+}\frac{\varphi_M(y_0)-\varphi_M(y_0+s\tilde y)}s
\\&=&-2\epsilon\beta r^2e^{-\beta r^2}\varphi_N(x_0)<0,
\end{eqnarray*}
completing the proof of \equ{eq ms3}.

The case for general $h$ can be done simply by replacing $h$ by $h^+$. In fact, since
    $w>0$ in $\Omega$, we have
$$(-\Delta)^\alpha_xw(x,y)+(-\Delta)_yw(x,y)+h^+(x,y)w(x,y)\ge0,\quad (x,y)\in \Omega$$
and similarly we obtain that
$$(-\Delta)^\alpha_xv(x,y)+(-\Delta)_yv(x,y)+h^+(x,y)v(x,y)\le0,\quad (x,y)\in D,$$
so we may proceed as before to get \equ{eq ms3} and the  proof is complete.
\hfill$\Box$\\

In order to state   the Strong Maximum Principle to be used in our
moving planes procedure, it is convenient to consider  property $(P)$:
\begin{itemize}
\item[$(P)$] We say that a
  function
$w:\tilde{\Omega}\to \R$ satisfies property $(P)$ if whenever  $(x_0,y_0)\in \Omega$ such that
 $$0=w(x_0,y_0)=\inf_{(x,y)\in\Omega}w(x,y),$$ then
$$w(x,y_0)\equiv0,\quad  \forall x\in\R^{N}.$$
\end{itemize}
The following lemma is in preparation of the strong maximum principle.
\begin{lemma}\label{proposition 2.1}
Let $\Omega$ be an open set in $\R^N\times\R^M$
and $w$ have property $(P)$. We denote
\begin{equation}\label{eq yx1}
\Omega_0=\{(x,y)\in\Omega: w(x,y)=\inf_\Omega w=0\}.
\end{equation}
If
$\O\not=\Omega_0\subsetneqq\Omega$, then $\Omega\setminus\Omega_0$
satisfies interior cylinder condition at any point
$(x_0,y_0)\in\partial\Omega_0\cap\Omega$.
\end{lemma}
{\bf Proof.} Since  $\O\not=\Omega_0\subsetneqq\Omega$,
we have that $\O\not=\partial\Omega_0\cap\Omega\subset\partial(\Omega\setminus\Omega_0)$. For any
$(x_0,y_0)\in\partial\Omega_0\cap\Omega$, let us denote $r=\frac{1}{4}dist((x_0,y_0),\partial\Omega)$ and let
$\tilde{y}\in\R^M$ such that $(x_0,\tilde{y})\in\Omega\setminus\Omega_0$
and $|\tilde{y}-y_0|=r$.
Since $w$ has property $(P)$, then $w=0$ in $\tilde \Omega_0$,
where $\tilde \Omega_0$ is the extension of $\Omega_0$ in $x$-direction and as
$\Omega\setminus\Omega_0$ is open, we have that   $B_r^N(x_0)\times B_r^M(\tilde y)\subset\Omega\setminus\Omega_0$.
Therefore, $\Omega\setminus\Omega_0$
satisfies interior cylinder condition at
$(x_0,y_0)\in\partial\Omega_0\cap\Omega$. \hfill$\Box$

\begin{teo}\label{teo SMP} [Strong Maximum Principle]\quad
Let $\Omega$ be an open set of $\R^N\times\R^M$, the function $h\in L_{loc}^\infty(\Omega)$ and $w\in C(\bar \Omega)\cap L^\infty(\tilde \Omega)$  has the
property $(P)$ satisfying
\begin{equation}\label{eq ms8}
\mathcal{L}w\ge0\quad {\rm{in}}\  \Omega \quad{\rm{and}}\quad
w\ge0\quad{\rm{in}}\ \Omega.
\end{equation}
Assume that
$\Omega_0\not=\O$ defined by (\ref{eq yx1}) and there exists some
$(x_0,y_0)\in\partial\Omega_0\cap \Omega$ such that (\ref{eq w2}) holds in
corresponding $D$.

 Then $w$ must be 0 in $\tilde \Omega$.
\end{teo}

\noindent{\bf Proof.} Assume that
$\Omega_0\not=\Omega.$ By Lemma \ref{proposition 2.1}, $\Omega\setminus\Omega_0$ satisfies interior cylinder condition at
$(x_0,y_0)\in\partial\Omega_0\cap\Omega$ and then $w(x_0,y_0)=0$ by $w\in C(\bar \Omega)$ and the definition of $\Omega_0$.
Furthermore,  we observe that $\bar D$ is compact in $\Omega$ and then $h\in L^\infty(\bar D)$.
 Using Lemma \ref{hopf}, we obtain  \equ{eq ms3},   which is impossible by the fact of $ w(x_0,y_0)=\inf_\Omega w=0$.
Therefore, $\Omega_0=\Omega$, i.e. $w\equiv0 $ in $\Omega$.  Since $w$ has property $(P)$, then $w\equiv0$ in $\tilde \Omega$.
 \hfill$\Box$

 \setcounter{equation}{0}
\section{Decay estimate}
\subsection{Proof of Theorem \ref{teo r1}}

In this subsection, we prove Theorem \ref{teo r1} on  decay estimates for positive classical solutions of equation (\ref{eq wyq1}).
The main work is to construct appropriate super and sub solutions and then the decay estimate is derived by
 Lemma \ref{teo MP}.

Before proving Theorem \ref{teo r1}, we introduce some computations gathered in the next proposition.
For $\alpha\in(0,1)$ and $\mu>0$, we define the function $ \psi_\mu:\R^N\to \R$ as follows
 \begin{equation}\label{eq2.1}
 \psi_\mu(x)=
\left\{ \arraycolsep=1pt
\begin{array}{lll}
 \mu^{-N-2\alpha}, \quad \ \  |x|<\mu,\\[2mm]
|x|^{-N-2\alpha},\quad \  |x|\geq\mu.
\end{array}
\right.
\end{equation}

\begin{proposition}\label{pr 2.1}
For any $\mu>0$, there exists $R_0>3\mu$ and $c>0$, independent of $\mu$, such that
\begin{equation}\label{3.07}
-c\mu^{-2\alpha}\psi_\mu(x)\le (-\Delta)^\alpha \psi_\mu(x)\le -c^{-1}\mu^{-2\alpha}\psi_\mu(x),\quad x\in B_{R_0}^c.
\end{equation}
\end{proposition}
{\bf Proof.} We consider along the proof that $\mu>0$ and $x\in \R^N $ satisfies  $|x|>3\mu$. We define
\begin{eqnarray*}
A(\mu,x,z)=\frac{  \psi_\mu(x+z)+ \psi_\mu(x-z)-2 \psi_\mu(x)}{|z|^{N+2\alpha}},\quad\quad z\in\R^N
\end{eqnarray*}
and we observe that
\begin{eqnarray}
(-\Delta)^\alpha \psi_\mu(x)=-\frac{1}{2}\int_{\R^N} A(\mu,x,z)dz.\label{eq 1.100}
\end{eqnarray}
Now we compute the integral above by decomposing the domain in various pieces.  First we consider the integral over $B_{\frac{|x|}{3}}(0)$. We observe that  $|x\pm z|\geq \mu$ for all  $z\in B_{\frac{|x|}{3}}(0)$, then by (\ref{eq2.1}) we obtain
\begin{eqnarray}
 |\int_{B_{\frac{|x|}{3}}(0)} A(\mu,x,z)dz|\nonumber
&=&|\int_{B_{\frac{|x|}{3}}(0)}\frac{|x+z|^{-N-2\alpha}+|x-z|^{-N-2\alpha}-2|x|^{-N-2\alpha}}{|z|^{N+2\alpha}}dz|\nonumber
\\&=&|x|^{-N-4\alpha}|\int_{B_{\frac{1}{3}}(0)}\frac{|z+e_x|^{-N-2\alpha}+|z-e_x|^{-N-2\alpha}-2}{|z|^{N+2\alpha}}dz|\nonumber
\\&\leq& c_1|x|^{-N-4\alpha}\int_{B_{\frac{1}{3}}(0)}\frac{|z|^2}{|z|^{N+2\alpha}}dz
\leq c_2|x|^{-N-4\alpha},\label{3.01}
\end{eqnarray}
where $e_x=\frac x{|x|}$ and $c_1,c_2>0$ are independent of $\mu$.
Next we consider the integral over $ B_{\frac{|x|}{3}}(x)\setminus B_\mu(x)$.
We observe that for all $z\in B_{\frac{|x|}{3}}(x)\setminus B_\mu(x)$ we have  $|x+z|\geq |x-z|\geq \mu$ and then we obtain
\begin{eqnarray*}
&& \int_{B_{\frac{|x|}{3}}(x)\setminus B_\mu(x)} A(\mu,x,z)dz
\\&=&\int_{B_{\frac{|x|}{3}}(x)\setminus B_\mu(x)} \frac{|x+z|^{-N-2\alpha}+|x-z|^{-N-2\alpha}-2|x|^{-N-2\alpha}}{|z|^{N+2\alpha}}dz
\\&=&|x|^{-N-4\alpha}\int_{B_{\frac{1}{3}}(e_x)\setminus B_{\frac\mu{|x|}}(e_x)}\frac{ |z+e_x|^{-N-2\alpha}+|z-e_x|^{-N-2\alpha}-2}{|z|^{N+2\alpha}}dz
\\&\leq& c_3|x|^{-N-4\alpha}\int_{B_{\frac{1}{3}}(e_x)\setminus B_{\frac\mu{|x|}}(e_x)}{|z-e_x|^{-N-2\alpha}}dz
\leq c_4\mu^{-2\alpha}|x|^{-N-2\alpha},
\end{eqnarray*}
where  the first inequality holds since  $|z+e_x|\geq |z-e_x|$ for $z\in B_{\frac{1}{3}}(e_x)\setminus B_{\frac\mu{|x|}}(e_x)$
and  $|z|\geq\frac{2}{3}$ for $z\in B_{\frac{1}{3}}(e_x)$.
For the inequality on the  other side, we obtain
\begin{eqnarray*}
 &&\int_{B_{\frac{|x|}{3}}(x)\setminus B_\mu(x)} A(\mu,x,z)dz
\\&=&|x|^{-N-4\alpha}\int_{B_{\frac{1}{3}}(e_x)\setminus B_{\frac\mu{|x|}}(e_x)}\frac{|z+e_x|^{-N-2\alpha}+|z-e_x|^{-N-2\alpha}-2}{|z|^{N+2\alpha}}dz
\\&\geq& |x|^{-N-4\alpha}(\int_{B_{\frac{1}{3}}(e_x)\setminus B_{\frac\mu{|x|}}(e_x)}\frac{|z-e_x|^{-N-2\alpha}}{|z|^{N+2\alpha}}dz
-\int_{B_{\frac{1}{3}}(e_x)}\frac{2}{|z|^{N+2\alpha}}dz)
\\&\geq &c_5|x|^{-N-4\alpha}\int_{B_{\frac{1}{3}}(e_x)\setminus B_{\frac\mu{|x|}}(e_x)}|z-e_x|^{-N-2\alpha}dz-c_{6}|x|^{-N-4\alpha}
\\&\geq& c_{7}\mu^{-2\alpha}|x|^{-N-2\alpha}- c_{8}|x|^{-N-4\alpha},
\end{eqnarray*}
where the second inequality holds by $|z|\leq\frac{4}{3}$ for $z\in B_{\frac{1}{3}}(e_x)$. Consequently,
\begin{equation}\label{3.030}
c_{7}\mu^{-2\alpha}|x|^{-N-2\alpha}- c_{8}|x|^{-N-4\alpha}\leq \int_{B_{\frac{|x|}{3}}(x)\setminus B_\mu(x)} A(\mu,x,z)dz
\leq c_{4}\mu^{-2\alpha}|x|^{-N-2\alpha},
\end{equation}
where the constants $c_4,c_{7},c_{8}>0$ are independent of $\mu$.
The estimate for the integral over $B_{\frac{|x|}{3}}(-x)\setminus B_\mu(-x)$ is similar.

Next we consider the integral over $ B_\mu(x)$. We observe  that, for $z\in  B_\mu(x)$ we have
since $|x+z|>\mu>|x-z|$ and $|z|\geq |x|-\mu\geq\frac{2|x|}{3}$, thus
\begin{eqnarray}
\nonumber &&\int_{ B_\mu(x)} A(\mu,x,z)dz
=\int_{ B_\mu(x)}\frac{|x+z|^{-N-2\alpha}+\mu^{-N-2\alpha}-2|x|^{-N-2\alpha}}{|z|^{N+2\alpha}}dz\nonumber
\\&&\le2\int_{B_\mu(x)}\frac{\mu^{-N-2\alpha}}{|z|^{N+2\alpha}}dz\nonumber
\leq  c_{9}\mu^{-2\alpha}(|x|-\mu)^{-N-2\alpha}
\leq c_{10}\mu^{-2\alpha}|x|^{-N-2\alpha}\label{3.04}\nonumber
\end{eqnarray}
and, for the other inequality
\begin{eqnarray}
\nonumber \int_{ B_\mu(x)} A(\mu,x,z)dz
&\ge&\int_{B_\mu(x)}\frac{-2|x|^{-N-2\alpha}}{|z|^{N+2\alpha}}dz
\\&\geq& -c_{11}\mu^{N}|x|^{-N-2\alpha}(|x|-\mu)^{-N-2\alpha}\nonumber
\geq - c_{12}|x|^{-N-4\alpha},\label{3.05}\nonumber
\end{eqnarray}
where  $c_{9}, c_{10}, c_{11}$ and $c_{12}$ are positive constant independent of $\mu$.
Therefore,
\begin{eqnarray}
 - c_{12}|x|^{-N-4\alpha}\leq \int_{ B_\mu(x)} A(\mu,x,z)dz
\leq  c_{10}\mu^{-2\alpha}|x|^{-N-2\alpha}.\label{3.050}
\end{eqnarray}
The integral over  $ B_\mu(-x) $ is exactly the same.
Finally, we consider the complementary integral over  $D(x)=\R^N\setminus(B_{\frac{|x|}{3}}(0)\cup {B_{\frac{|x|}{3}}(x)}\cup {B_{\frac{|x|}{3}}(-x)})$.
For $|x|>3\mu$ and $z\in D(x)$, we have that $|x\pm z|\geq \frac{|x|}{3}$, thus
\begin{eqnarray}
 |\int_{D(x)} A(\mu,x,z)dz|
&\leq& \int_{D(x)}\frac{ |x+z|^{-N-2\alpha}+|x-z|^{-N-2\alpha}+2|x|^{-N-2\alpha}}{|z|^{N+2\alpha}}dz\nonumber
\\&\leq& c_{13}|x|^{-N-2\alpha}\int_{ \R^N\setminus B_{\frac{|x|}{3}}(0) }\frac{1}{|z|^{N+2\alpha}}dz	\nonumber
\\&\leq& c_{14}|x|^{-N-4\alpha},\label{3.06}
\end{eqnarray}
where $c_{13}>0$ and $c_{14}>0$ are independent of $\mu$.
Therefore, by (\ref{3.01})-(\ref{3.06}), there exist $c_{15},c_{16}>1$ independent of $\mu$ such that
\begin{eqnarray*}
&& c_{15}^{-1}\mu^{-2\alpha}|x|^{-N-2\alpha}- c_{15}|x|^{-N-4\alpha}\leq \int_{\R^N} A(\mu,x,z)dz
\\&&\leq c_{16}\mu^{-2\alpha}|x|^{-N-2\alpha}+c_{16}|x|^{-N-4\alpha}
\leq  c_{15}\mu^{-2\alpha}|x|^{-N-2\alpha},
\end{eqnarray*}
where we used that  $|x|>3\mu$. Choosing $R_0>3\mu$ such that
$c_{15}^{-1}\mu^{-2\alpha}- c_{15}|x|^{-2\alpha}\geq\frac{ 1}{2}c_{15}^{-1}\mu^{-2\alpha}$
for $|x|\geq R_0$,
together with (\ref{eq 1.100}), we obtain   (\ref{3.07}).
\hfill$\Box$

\vspace{3mm}

In what follows we provide a proof of our first theorem on the decay of the positive solutions of our equation.

\noindent{\bf Proof of Theorem \ref{teo r1}.}
 By definition of $A$ and $B$ in \equ{eq 011j}, for
 any $\epsilon>0$,
there exits
$\delta_\epsilon>0$ such that
\begin{equation}\label{eq11}
(B-\epsilon^2)t\leq f(t)\leq (A+\epsilon^2)t,\quad \ \forall \ t\in(0,\delta_\epsilon).
\end{equation}
 Since $u$ is a positive solution of (\ref{eq wyq1}) vanishing  at infinity,
 there exists $R_\epsilon>0$ such that $0<u(x,y)<\delta_\epsilon$ for any $(x,y)\in B^c_{R_\epsilon}$. Therefore,
\begin{equation}\label{f2}
(-\Delta)_x^\alpha u+(-\Delta)_yu+ (1-A-\epsilon^2)u\leq 0 \quad {\rm{in}}\ B^c_{R_\epsilon}
\end{equation}
and
\begin{equation}\label{f02}
(-\Delta)_x^\alpha u+(-\Delta)_yu+(1-B+\epsilon^2)u\geq0 \quad {\rm{in}}\ B^c_{R_\epsilon}.
\end{equation}
Next we define the function $\phi_\nu:\R^M\to\R$ as
$\phi_\nu(y)=e^{-\nu|y|},$
where $\nu>0$ and we find  that
for $ y\in\R^M\setminus\{0\}$,
\begin{equation}\label{eq2.02}
-\Delta \phi_\nu(y)=\nu\left(\frac{M-1}{|y|}-\nu\right)\phi_\nu(y).
 \end{equation}

\noindent\textbf{ Step 1.}\emph{ There exists $C(\epsilon)>1$ such that
\begin{equation}\label{eq mt1}
u(x,y)\leq C(\epsilon) e^{-\theta_1|y|},\quad \quad (x,y)\in\R^N\times\R^M.
 \end{equation}
}
To prove \equ{eq mt1}  we let
$U_1(x,y)=\phi_{\theta_1}(y)$, for  $(x,y)\in\R^N\times\R^M$
and then,  by (\ref{eq2.02}), we have
\begin{eqnarray}
&&(-\Delta)_x^\alpha U_1+(-\Delta)_yU_1+(1-A-\epsilon^2)U_1
\nonumber\\&&
=\left[\theta_1\left(\frac{M-1}{|y|}-\theta_1\right)+1-A-\epsilon^2\right]U_1\ge 0,\label{eq r212}
\end{eqnarray}
if $\epsilon\le \sqrt{1-A} $. By definition of $U_1$ and $\phi_{\theta_1}$ we have that
$U_1=1$ in $\R^N\times\{0\}$ and $U_1\geq e^{-\theta_1 R_\epsilon}$ in $\bar B_{R_\epsilon}$ and, since
 $u$ is bounded,
 there exists
$\rho_1>0$ depending on $\epsilon$, such that
$$W_1=\rho_1 U_1-u\ge 0\quad {\rm{in}}\ \ \bar B_{R_\epsilon}\cup(\R^N\times\{0\}).$$
Combining (\ref{f2})  with (\ref{eq r212}), we obtain
$$(-\Delta)_x^\alpha W_1+(-\Delta)_yW_1+(1-A-\epsilon^2)W_1\ge0\quad {\rm{in}}\ \ \bar B_{R_\epsilon}^c\cap(\R^N\times\{0\})^c.$$
By Lemma \ref{teo MP}, this implies that
$W_1\ge0$ in $\R^N\times\R^M$
and then
\begin{equation}\label{eq mt112}
u(x,y)\leq \rho_1 U_1(x,y)=\rho_1 \phi_{\theta_1}(y)=\rho_1 e^{-\theta_1|y|},\quad  (x,y)\in\R^N\times\R^M. \end{equation}

\noindent\textbf{ Step 2.} \emph{There exists $C(\epsilon)>1$ such that
\begin{equation}\label{eq mt2}
u(x,y)\leq C(\epsilon) |x|^{-N-2\alpha}, \quad \quad (x,y)\in\R^N\times\R^M.
 \end{equation}}
Let $c$ and $R_0$ be as in Proposition \ref{pr 2.1}
$\mu=({c}/({2\epsilon\sqrt{(1-A)}-2\epsilon^2}))^{\frac1{2\alpha}}$
and consider the function
$
U_2(x,y)=\psi_\mu(x),$ for $(x,y)\in\R^N\times\R^M.$
Then, by (\ref{3.07}), we have for all $(x,y)\in (B^N_{R_0}(0))^c\times\R^M$ that
\begin{eqnarray}
&&
(-\Delta)_x^\alpha U_2+(-\Delta)_y U_2+(1-A-\epsilon^2)U_2
\nonumber\\&&
\ge(-c\mu^{-2\alpha}+1-A-\epsilon^2)U_2
\ge0\label{eq r2220}
\end{eqnarray}
for  $0<\epsilon<\sqrt{1-A}$. Let us
denote $W_2=\rho_2 U_2-u$, where $\rho_2>0$ is such that
$$W_2\ge \rho_2(R_0+R_\epsilon)^{-N-2\alpha}-u\ge 0\quad {\rm{in}}\ \bar B_{R_\epsilon}\cup ( \overline{B_{R_0}^N(0)}\times \R^M).$$
Combining (\ref{f2})  with (\ref{eq r2220}), we obtain that
$$\!(-\Delta)_x^\alpha W_2+(-\Delta)_yW_2+(1-A-\epsilon^2)W_2\ge0\quad {\rm{in}}\ \ \bar B_{R_\epsilon}^c\cap ( \overline{B_{R_0}^N(0)}\times \R^M)^c.$$
By Lemma \ref{teo MP}, we have that
$W_2=\rho_2 U_2-u\ge0$ {{in}} $R^N\times\R^M$ and then,
 for all $(x,y)\in\R^N\times\R^M$,
$$u(x,y)\leq \rho_2 U_2(x,y)=\rho_2 \psi_\mu(x)\leq \rho_2|x|^{-N-2\alpha}.$$

\noindent\textbf{Step 3.} \emph{There exists  $C(\epsilon)>1$ such that
\begin{equation}\label{eq mt8}
u(x,y)\leq C(\epsilon) |x|^{-N-2\alpha}e^{-\theta_1|y|}, \quad\ \ (x,y)\in\R^N\times\R^M.
 \end{equation}
}
Let us  consider the function
$V(x,y)=\psi_\mu(x) \phi_{\theta_1}(y),$ for $(x,y)\in\R^N\times\R^M,$
with $\mu$ as defined above. From (\ref{3.07}) and (\ref{eq2.02}), we have that
\begin{eqnarray}
&&(-\Delta)_x^\alpha V+(-\Delta)_yV+(1-A-\epsilon^2)V\nonumber
\\&&\ge\left[-c\mu^{-2\alpha}+\theta_1\left(\frac{M-1}{|y|}-\theta_1\right)+1-A-\epsilon^2\right]V
\ge0,\label{eq r2220m}
\end{eqnarray}
for $(x,y)\in (B^N_{R_0}(0))^c\times(\R^M\setminus\{0\})$ and assuming that $0<\epsilon<\sqrt{1-A}$.
Since $u$, $V$ are bounded in $\bar B_{R_\epsilon}$ and $V$ is positive, there is $\bar\rho_1>0$ large such that
$$
\bar\rho_1 V-u\geq0\quad {\rm{in}}\ \ \bar B_{R_\epsilon}.
$$
By \equ{eq mt1} and (\ref{eq mt112}), we may choose $\bar\rho_2>0$  such that
\begin{eqnarray*}
\bar\rho_2 V-u &\geq &\bar\rho_2R_0^{-N-2\alpha}\phi_{\theta_1}(y)-u\geq0\quad {\rm{in}}\ \ \overline{B^N_{R_0}(0)}\times\R^M\quad \mbox{and}\\
\bar\rho_2 V-u & \geq &\bar\rho_2\psi_\mu(x) -u\geq0\quad {\rm{in}}\ \ \R^N\times\{0\}.
\end{eqnarray*}
Taking $\bar\rho=\max\{ \bar\rho_1, \bar\rho_2\}$, defining $W=\bar\rho V-u$ and combining (\ref{f2})  with (\ref{eq r2220m}), we have that
$$
W \geq0  \ \ {\rm{in}}\ \  \bar B_{R_\epsilon}\cup(\overline{B^N_{R_0}(0)}\times\R^M)\cup(\R^N\times\{0\})\quad \mbox{and}
$$
$$(-\Delta)_x^\alpha W+(-\Delta)_yW+(1-A-\epsilon^2)W \ge0 \ \  {\rm{in}}\ \ \bar B_{R_\epsilon}^c\cap ((B^N_{R_0}(0))^c\times(\R^M\setminus\{0\})).
$$
Then, by Lemma \ref{teo MP}, we have that
$\bar\rho V-u\ge0$ in $ \R^N\times\R^M.$
Thus,  there exists  $ C(\epsilon)>1$  such that
$$
u(x,y)\leq C(\epsilon)\psi_\mu(x) \phi_{\theta_1}(y)\leq   C(\epsilon) |x|^{-N-2\alpha}e^{-\theta_1|y|},\qquad (x,y)\in\R^N\times\R^M.
$$

\noindent\textbf{ Step 4.}  \emph{There exists $C_1(\epsilon)>0$ and $R>0$ such that
\begin{equation}\label{eq mt113}
u(x,y)\geq C_1(\epsilon) e^{-\theta_2|y|},\quad \quad (x,y)\in \overline{B_{R}^N(0)}\times \R^M.
 \end{equation}
}
Let $R_0$ be as in Proposition \ref{pr 2.1} and let $R>
R_0$  such that $\lambda_1<\epsilon^2$, where  $\lambda_1$ is the first eigenvalue of the fractional Dirichlet problem (\ref{eq w4}) with $x_0=0$ and $r=4R$.
 Let  $\varphi_N$ be the first eigenfunction of (\ref{eq w4}) and define
$V_1(x,y)=\varphi_N(x)\phi_{\theta_2}(y)$ for $(x,y)\in\R^N\times\R^M.$
From  (\ref{eq w4}) and (\ref{eq2.02}), for  $(x,y)\in B_{2R}^N(0)\times (B_{R_1}^M(0))^c$ with $R_1=\frac{M-1}{\epsilon}$,    we have
\begin{eqnarray}
(-\Delta)_x^\alpha V_1+(-\Delta)_yV_1+(1-B+\epsilon^2)V_1\nonumber \\
=\left[\lambda_1+{\theta_2}\left(\frac{M-1}{|y|}-{\theta_2}\right)+1-B+\epsilon^2\right]V_1\nonumber \\
\leq[\epsilon^2+{\theta_2}(\epsilon-{\theta_2})+1-B+\epsilon^2]V_1\leq0,\label{eq r221}
\end{eqnarray}
if $\epsilon<\sqrt{1-B}$.
Let us define $w_1=u-r_1 V_1$, where $r_1>0$ is
such that
$$w_1\ge 0\quad {\rm{in}}\ \overline B_{R_\epsilon}\cup(\overline{B_{2R}^N(0)\times B_{R_1}^M(0)})$$
and observe that  $w_1\ge 0$ in $(B_{2R}^N(0))^c\times \R^M$ since $V_1=0$.
Combining (\ref{f02})  with (\ref{eq r221}), we obtain that
$$(-\Delta)_x^\alpha w_1+(-\Delta)_yw_1+(1-B+\epsilon^2)w_1\ge0\quad {\rm{in}}\ \ (B_{2R}^N(0)\times (B_{R_1}^M(0))^c)\cap B_{R_\epsilon}^c$$
and then, by Lemma \ref{teo MP}, we have that
$$w_1=u-r_1 V_1\ge0\quad {\rm{in}}\ \ \R^N\times\R^M.$$
Since $\varphi_N$ is classical solution of (\ref{eq w4}) with $r=4R$ and $x_0=0$ then  $\varphi_N(x)$ is positive in
 $\overline{B_{R}^N(0)}\subset\R^N$, we can finally choose
 $C_1(\epsilon)>0$ such that
\begin{equation}\label{eq 1.12}
u(x,y)\geq r_1\varphi_N(x) \phi_{\theta_2}(y)\geq C_1(\epsilon)e^{-\theta_2|y|},\ \  \forall  (x,y)\in \overline{B_{R}^N(0)}\times \R^M.
 \end{equation}

\noindent\textbf{Step 5.}\emph{ There exists  $ C_1(\epsilon)>0$ such that, for $R$ and $R_1$ as in Step 4,
\begin{equation}\label{eq 1.13}
u(x,y)\geq C_1 |x|^{-N-2\alpha}, \quad\quad (x,y)\in  (B_{R}^N(0))^c\times \overline{B_{R_1}^M(0)}.
 \end{equation}
}
 To prove this,
we define
$V_2(x,y)=\psi_\mu(x)\eta_M(y)$ for $(x,y)\in\R^N\times\R^M,
$
where $\eta_M$ is the solution of
\begin{equation}\label{eq w42}
\left\{ \arraycolsep=1pt
\begin{array}{lll}
-\Delta \eta_M(y)=\bar\lambda_1 \eta_M(y) ,\ \ \ \ &
y\in B_{R_2}^M(0),\\[2mm]
\eta_M(y)=0,& y\in(B_{R_2}^M(0))^c,
\end{array}
\right.
\end{equation}
with $R_2>R_1$  such that $\bar\lambda_1 <\epsilon^2$. Here $\mu=[c(1-B+2\epsilon^2)]^{\frac{-1}{2\alpha}}$ with  $c$ as in Proposition \ref{pr 2.1} and $\psi_\mu$ defined in (\ref{eq2.1}).
By (\ref{3.07}) and (\ref{eq w42}), for $(x,y)\in ((B_{R}^N(0))^c\times\R^M)\cap (\R^N\times B_{R_2}^M(0))$, we have that
\begin{eqnarray}
&&(-\Delta)_x^\alpha V_2+(-\Delta)_yV_2+(1-B+\epsilon^2)V_2\nonumber
\\&&\le(-c^{-1}\mu^{-2\alpha}+\bar\lambda_1+1-B+\epsilon^2)V_2
=0.\label{eq r224}
\end{eqnarray}
Let $w_2=u-r_2 V_2$, with $r_2 >0$  such that
$$w_2\ge 0\quad {\rm{in}}\ \bar  B_{R_\epsilon}\cup\overline{(B_{R}^N(0)}\times\R^M)\cup (\R^N\times (B_{R_2}^M(0))^c).$$
Combining (\ref{f02})  with (\ref{eq r224}), we obtain that
$$(-\Delta)_x^\alpha w_2+(-\Delta)_yw_2+(1-B+\epsilon^2)w_2\ge0$$
in $B_{R_\epsilon}^c\cap((B_{R}^N(0))^c\times\R^M)\cap (\R^N\times B_{R_2}^M(0)).$
By Lemma \ref{teo MP}, we have then
$$w_2=u-r_2 V_2\ge0\quad {\rm{in}}\ \ \R^N\times\R^M.$$
Since $\eta_M$ is  positive
in $ \overline{B_{R_1}^M(0)}\subset{B_{R_2}^M(0)}$,  there exists  $C_1(\epsilon)>0$ such that for any $(x,y)\in (B_{R}^N(0))^c\times \overline{B_{R_1}^M(0)}$, we have that
$$
u(x,y)\geq r_2 \psi_\mu(x)\eta_M(y)\geq C_1(\epsilon)  |x|^{-N-2\alpha}.
$$

\noindent\textbf{ Step 6.} \emph{There exist $C_1(\epsilon)>0$  such that, for $R$ as in   Step 4,
\begin{equation}\label{eq mt8n}
u(x,y)\geq C_1(\epsilon)|x|^{-N-2\alpha} e^{-\theta_2|y|},\quad \quad (x,y)\in(B_{R}^N(0))^c\times\R^M.
 \end{equation}
 }
To prove this we let
$\tilde{V}(x,y)=\psi_\mu(x) \phi_{\theta_2}(y),$ for $(x,y)\in\R^N\times\R^M$
with $\mu$ as defined above.
Using (\ref{3.07}) and (\ref{eq2.02}), for $(x,y)\in (B^N_{R}(0))^c\times(B^M_{R_1}(0))^c$ with $R_1=\frac{M-1}{\epsilon}$, we have that
\begin{eqnarray}
&&(-\Delta)_x^\alpha \tilde{V}+(-\Delta)_y\tilde{V}+(1-B+\epsilon^2)\tilde{V}\nonumber
\\&&\le\left[-c^{-1}\mu^{-2\alpha}+{\theta_2}\left(\frac{M-1}{|y|}-{\theta_2}\right)+1-B+\epsilon^2\right]\tilde{V}\nonumber
\\&&\le[{\theta_2}(\epsilon-{\theta_2})+1-B+\epsilon^2]\tilde{V}
\le0,\label{eq r2220my}
\end{eqnarray}
if $0<\epsilon<\sqrt{1-B}$.
Since $u$ is positive  and  $V$ is bounded in $\overline B_{R_\epsilon}$,  we can choose $\tilde{r}_1>0$  such that
$$
u-\tilde{r}_1 V\geq0\quad {\rm{in}}\ \ \overline B_{R_\epsilon}.
$$
Since $\psi_\mu$ is bounded  in $\overline{B^N_{R}(0)}$, using
(\ref{eq 1.12}),  there exists  $\tilde{r}_2>0$ such that
$$
u-\tilde{r}_2 V\geq u-\tilde{r}_2c_1e^{-\theta_2|y|}\geq0\quad {\rm{in}}\ \ \overline{B^N_{R}(0)}\times\R^M,
$$
and by (\ref{eq 1.13}),  there exists $\tilde{r}_3>0$  such that
$$
u-\tilde{r}_3 V\geq u-\tilde{r}_3|x|^{-N-2\alpha}\geq0\quad {\rm{in}}\ \ (B_{R}^N(0))^c\times \overline{B_{R_1}^M(0)}.
$$
Taking $\tilde{r}=\min\{ \tilde{r}_1, \tilde{r}_2, \tilde{r}_3\}$ and combining (\ref{f02})  with (\ref{eq r2220my}), we obtain that
\begin{eqnarray*}
w=u-\tilde{r} V\geq0\quad {\rm{in}}\ \ \bar B_{R_\epsilon}\cup(\overline{B^N_{R}(0)}\times\R^M)\cup((B_{R}^N(0))^c\times \overline{B_{R_1}^M(0)})\ \ \mbox{and}\\
(-\Delta)_x^\alpha w+(-\Delta)_yw+(1-B+\epsilon)w\ge0\quad {\rm{in}}\ \ \bar B_{R_\epsilon}^c\cap ((B^N_{R}(0))^c\times(B^M_{R_1}(0))^c).
\end{eqnarray*}
Thus Lemma \ref{teo MP}, we have that
$w\ge0$ in $\R^N\times\R^M$
and then (\ref{eq mt8n}) holds.

Finally,  Step 1 $-$  Step 6 completes  the proof. \hfill$\Box$

\subsection{Proof of Theorem \ref{teo jxcwr1}}

This subsection is devoted to prove Theorem \ref{teo jxcwr1}. Our proof is based on the fundamental solution of the   mixed integro-differential operator.
We first study the
fundamental solution $\mathcal{K}$ for
$$(-\Delta)_x^{\alpha} u+(-\Delta)_y u+ u=0\quad {\rm in}\quad \R^N\times(\R^M\setminus\{0\}),$$
which can be characterized by
\begin{equation}\label{eq hk2}
\mathcal{K}(x,y)=\int_0^\infty e^{-t}\mathcal{H}(x,y,t)dt,
\end{equation}
where
 \begin{equation}\label{eq hk1}
\mathcal{H}(x,y,t)=\int_{\R^M}\int_{\R^N}e^{-2\pi i(x,y)\cdot(\xi_1,\xi_2)-t(|\xi_1|^{2\alpha}+|\xi_2|^2)}d\xi_1d\xi_2.
\end{equation}
In fact, for $\phi\in\mathcal{S}$, we have that
$$
 \arraycolsep=1pt
\begin{array}{lll}
\langle\mathcal{K},\phi\rangle=\int_{\R^{N+M}}\int_0^\infty\int_{\R^{N+M}} e^{-2\pi i(x,y)\cdot(\xi_1,\xi_2)-t(|\xi_1|^{2\alpha}+|\xi_2|^2+1)}\phi(x,y)d\xi_1d\xi_2dtdxdy
\\[2mm]\phantom{-}
=\int_{\R^{N+M}}\left[\int_0^\infty e^{-t(|\xi_1|^{2\alpha}+|\xi_2|^2+1)}dt\int_{\R^{N+M}} e^{-2\pi i(x,y)\cdot(\xi_1,\xi_2)}\phi(x,y)dxdy
\right]d\xi_1d\xi_2
\\[2mm]\phantom{-}
=\int_{\R^{N+M}}\left[\frac{1}{|\xi_1|^{2\alpha}+|\xi_2|^2+1}\int_{\R^{N+M}} e^{-2\pi i(x,y)\cdot(\xi_1,\xi_2)}\phi(x,y)dxdy
\right]d\xi_1d\xi_2
\\[2mm]\phantom{-}
=\left\langle\frac{1}{|\xi_1|^{2\alpha}+|\xi_2|^2+1},\mathcal{F}\phi\right\rangle.
\end{array}
$$
Next we want to find some properties of  $\mathcal{H}$.  To this end, we consider
$$\mathcal{H}_\alpha(x,t)=\int_{\R^N}e^{-2\pi ix\cdot\xi_1-t|\xi_1|^{2\alpha}}d\xi_1  \quad
{\rm and} \quad \mathcal{H}_1(y,t)=\int_{\R^M}e^{-2\pi iy\cdot\xi_2-t|\xi_2|^2}d\xi_2.$$
It is well known that the function
 $\mathcal{H}_\alpha$  has the following  properties:
$$\mathcal{H}_\alpha(x,t)=t^{-\frac{N}{2\alpha}}\mathcal{H}_\alpha(t^{-\frac{1}{2\alpha}}x,1) \quad
{\rm and} \quad \lim_{|x|\to\infty}|x|^{N+2\alpha}\mathcal{H}_\alpha(x,1)=C,$$
where $C>0$,
which imply that there exists $c_1>0$ and $c_2>$ such that
\begin{equation}\label{eq hk30}
c_{1}\min\{t^{-\frac{N}{2\alpha}},t|x|^{-N-2\alpha}\}
\leq \mathcal{H}_\alpha(x,t)
\leq c_2\min\{t^{-\frac{N}{2\alpha}},t|x|^{-N-2\alpha}\},
\end{equation}
see \cite{REFERENCE-Fourier,FQT}.
By the definition of  $\mathcal{H}$,
we have that
\begin{equation}\label{eq hk3}
\mathcal{H}(x,y,t)=\mathcal{H}_\alpha(x,t)\mathcal{H}_1(y,t).
\end{equation}
Since we have
\begin{equation}\label{eq hk30jx1}
 \mathcal{H}_1(y,t)=(4\pi t)^{-\frac{M}{2}}e^{-\frac{|y|^2}{4t}},
 \end{equation}
see \cite{REFERENCE-Fourier}, together with (\ref{eq hk2})-(\ref{eq hk3}),  for $|y|>2$,
\begin{eqnarray*}
\mathcal{K}(x,y)&=&\int_0^\infty e^{-t}\mathcal{H}_\alpha(x,t)\mathcal{H}_1(y,t)dt
\\&\geq & c_{1}\int_0^\infty  e^{-t}\min\{t^{-\frac{N}{2\alpha}},t|x|^{-N-2\alpha}\}(4\pi t)^{-\frac M2}e^{-\frac{|y|^2}{4t}}dt
\\&\geq & c_{1}\int_{\frac{|y|}{2}}^{\frac{|y|}{2}+1}  e^{-t}\min\{t^{-\frac{N}{2\alpha}},t|x|^{-N-2\alpha}\}(4\pi t)^{-\frac M2}e^{-\frac{|y|^2}{4t}}dt
\\&\geq & c_{3}\min\{e^{-|y|}|y|^{-\frac{N}{2\alpha}-\frac M2},  |x|^{-N-2\alpha} e^{-|y|}|y|^{1-\frac{M}{2}} \},
\end{eqnarray*}
for some $c_{3}>0$.
On the other hand, since for  $n\geq 3$ we have
$$\int_0^\infty  e^{-t}(4\pi t)^{-\frac{n}{2}}e^{-\frac{|y|^2}{4t}}dt
\leq c_{4} e^{-|y|}|y|^{2-n}(1+|y|)^{\frac{n-3}{2}}$$
with $c_{4}>0$ (see \cite{REFERENCE-Fourier}), for $M\ge 5$ we have
 that
\begin{eqnarray*}
&&\mathcal{K}(x,y)=\int_0^\infty e^{-t}\mathcal{H}_\alpha(x,t)\mathcal{H}_1(y,t)dt
\\&\leq & c_{2}\int_0^\infty  e^{-t}\min\{t^{-\frac{N}{2\alpha}},t|x|^{-N-2\alpha}\}(4\pi t)^{-\frac M2}e^{-\frac{|y|^2}{4t}}dt
\\&\leq & c_{5}\min\{\int_0^\infty  e^{-t}(4\pi t)^{-\frac{N}{2\alpha}-\frac M2}e^{-\frac{|y|^2}{4t}}dt,|x|^{-N-2\alpha}\int_0^\infty  e^{-t}(4\pi t)^{1-\frac M2}e^{-\frac{|y|^2}{4t}}dt\}
\\&\leq & c_{6}\min\{ e^{-|y|}|y|^{2-\frac{N}{\alpha}-M}(1+|y|)^{\frac{N}{2\alpha}+\frac M2-\frac32},
|x|^{-N-2\alpha}e^{-|y|}|y|^{4-M}(1+|y|)^{\frac {M-5}2}\}
\end{eqnarray*}
Therefore, for $N\geq1$ and $M\geq5$, there exist $c_{8}>c_{7}>0$ such that
\begin{equation}\label{jx2234}
c_{7}\rho(x,y)\leq \mathcal{K}(x,y)\leq  c_{8}\rho(x,y)|y|^{\frac{1}{2}}, \quad (x,y)\in\R^N\times(B_2^M(0))^c,
\end{equation}
where $\rho(x,y)$ is defined in \equ{ro}. In what follows, we construct super and sub-solutions to
obtain the decay estimate given  in Theorem  \ref{teo jxcwr1}.

\medskip

\noindent{\bf Proof of Theorem \ref{teo jxcwr1}.}
By the estimate in Theorem \ref{teo r1}, we observe that, for constants   $c_{10}>c_{9}>0$  such that
$$
 c_{9}(1+|x|)^{-N-2\alpha}\leq u(x,y)\leq  c_{10}(1+|x|)^{-N-2\alpha}, \quad (x,y)\in\R^N\times{B_2^M(0)},
$$
so we only need to prove (\ref{jx220}) holds for $(x,y)\in\R^N\times{(B_2^M(0))^c}$.

\noindent {\bf Step 1}: {\it Lower bound.}
Let $\tilde{u}=\mathcal{K}*\chi_{B_1^N(0)\times {B_1^M(0)}}$, where $\chi_{B_1^N(0)\times B_1^M(0)}$ is the characteristic function of
 $B_1^N(0)\times B_1^M(0)$. By (\ref{jx2234}),   we have that
\begin{equation}\label{eq jx01y}
  \tilde{u}(x,y)\geq  c_{11}\min\{e^{-|y|}|y|^{-\frac{N}{2\alpha}-\frac M2},  (1+|x|)^{-N-2\alpha} e^{-|y|}|y|^{1-\frac{M}{2}} \},
\end{equation}
for all $(x,y)\in\R^N\times{(B_2^M(0))^c}$, where $c_{11}>0$.
By definition of $\tilde{u}$, we have
\begin{eqnarray*}
(-\Delta)_x^\alpha \tilde{u}+(-\Delta)_y \tilde{u}+\tilde{u}=0\quad {\rm in}\ \  \R^N\times(\R^M\setminus\{0\})\setminus (B_1^N(0)\times B_1^M(0))
\end{eqnarray*}
and, by  (\ref{jx2234}) and  Theorem \ref{teo r1},  there exists $c_{12}>0$ such that
$u \geq c_{11}\tilde{u}$ in $\R^N\times\{y\in\R^M:|y|=2\}$.
Since $f$ is nonnegative, we use the Comparison Principle to obtain that, for any $(x,y)\in \R^N\times(B_2^M(0))^c$
 \begin{eqnarray*}
u(x,y) \geq  c_{11}\tilde{u}(x,y)
\geq c_{12}\min\{e^{-|y|}|y|^{-\frac{N}{2\alpha}-\frac M2},  (1+|x|)^{-N-2\alpha} e^{-|y|}|y|^{1-\frac{M}{2}} \}.
\end{eqnarray*}

\noindent
{\bf Step 2:} {\it Upper bound.} For $y\in \R^M$ with $|y|\geq2$, there exists $1\le i \le M$ such that $|y_i|>1$,
we may assume that $y_1>1$. Let $\bar u(x,y)=\mathcal{K}(x,y)(1-|y_1|^{-1})$, then by direct computation
 \begin{eqnarray*}
(-\Delta)_y\bar u&=& (1-|y_1|^{-1}) (-\Delta)_y\mathcal{K}-2y_1^{-2}\partial_{y_1}\mathcal{K}+2\mathcal{K}y_1^{-3}
\\&\geq & (-\Delta)_y\mathcal{K}(1-|y_1|^{-1})+2\mathcal{K}y_1^{-3},
\end{eqnarray*}
where the last inequality holds since  $y_1>0$ and
$\partial_{y_1}\mathcal{K}
<0.$
Therefore, by (\ref{jx2234}), we have that for $(x,y)\in\R^N\times(B_2^M(0))^c$,
 \begin{eqnarray}
&&(-\Delta)_x^\alpha \bar{u}(x,y)+(-\Delta)_y \bar{u}(x,y)+\bar{u}(x,y)\nonumber
\\&\geq&[(-\Delta)_x^\alpha \mathcal{K}+(-\Delta)_y  \mathcal{K}+ \mathcal{K}](1-|y_1|^{-1})+2\mathcal{K}(x,y)y_1^{-3}\nonumber
\geq  2\mathcal{K}(x,y)|y|^{-3}\nonumber
\\&\geq & 2c_{8}\min\{e^{-|y|}|y|^{-\frac{N}{2\alpha}-\frac M2-3},  |x|^{-N-2\alpha} e^{-|y|}|y|^{-\frac{M}{2}-2} \}. \label{eq 010}
\end{eqnarray}
Since $f(u)=O(u^p)$ near $u=0$ for some $p>1$, by  Theorem \ref{teo r1} with $\epsilon=\frac{p-1}{4p}$,
we have that
$$(-\Delta)_x^\alpha u+(-\Delta)_y u+u)=f(u)\leq c_{13}(1+|x|)^{-(N+2\alpha)p}e^{-\frac{3p+1}{4}|y|},$$
where $c_{13}>0$. We notice that $\frac{3p+1}{4}>1$. By  definition of $\bar u$,  (\ref{jx2234}) and Theorem \ref{teo r1} with
 $\epsilon=\frac{p-1}{4p}$, there exists $c_{14}>0$ such that
$u \leq c_{14}\bar u$ in $\R^N\times\{y\in\R^M:|y|=2\}$.
By Comparison Principle, we have that
 \begin{eqnarray*}
u(x,y)&\leq & c_{14}\bar u(x,y)\leq c_{14}\mathcal{K}(x,y)
\\&\leq & c_{15}\min\{e^{-|y|}|y|^{\frac{1}{2}-\frac{N}{2\alpha}-\frac M2},  (1+|x|)^{-N-2\alpha} e^{-|y|}|y|^{\frac{3}{2}-\frac{M}{2}} \}
\end{eqnarray*}
for all $(x,y)\in \R^N\times(B_2^M(0))^c$ and some $c_{15}>0$. This complete the proof.

 \hfill$\Box$

\setcounter{equation}{0}
\section{Symmetry results}

In this section,  we prove   Theorem
\ref{teo r1jt} by  moving planes method.
 Let $u$
be a classical positive solution of \equ{eq wyq1} and consider  first the $y$-direction. Let
$$\Sigma^{y_1}_\lambda=\{(x,y_1,y')\in \R^N\times\R\times\R^{M-1}\  |\  y_1>\lambda\},$$
$$T^{y_1}_\lambda=\{(x,y_1,y')\in  \R^N\times\R\times\R^{M-1}\  |\  y_1=\lambda\}$$
and $u_\lambda(x,y_1,y')=u(x,2\lambda-y_1,y')$ for $\lambda\in\R$.
We introduce  a preliminary inequality which plays a crucial role in
the procedure of moving planes.
\begin{lemma}\label{lemma dem4}
Under the assumptions of Theorem  \ref{teo r1jt}, for any
$\lambda\in\R$, there exists $c_{1}>0$, independent of
$\lambda$, such that
\begin{eqnarray*}
&&c_{1}(\int_{\Sigma^{y_1}_\lambda}|(u_\lambda-u)^+|^{\frac{2(N+M)}{N+M-2\alpha}}dxdy)^{\frac{N+M-2\alpha}{N+M}}
\\&\leq&
\int_{\Sigma^{y_1}_\lambda}[(-\Delta)_x^{\alpha}({u_\lambda}-u)+(-\Delta)_y({u_\lambda}-u)+({u_\lambda}-u)](u_\lambda-u)^+dxdy
<\infty.
\end{eqnarray*}
\end{lemma}
\noindent{\bf Proof.} {First we show that the integrals are finite.}
We observe that $u_\lambda$ satisfies the same equation (\ref{eq wyq1})
 as $u$ in $\Sigma^{y_1}_\lambda$. Taking $(u_\lambda-u)^+$ as test function in the equations for
$u$ and $u_\lambda$,  subtracting and integrating in
$\Sigma^{y_1}_\lambda$, we find
\begin{eqnarray}
&&\int_{\Sigma^{y_1}_\lambda}[(-\Delta)_x^{\alpha}({u_\lambda}-u)+(-\Delta)_y({u_\lambda}-u)+({u_\lambda}-u)](u_\lambda-u)^+dxdy\nonumber
\\&&= \int_{\Sigma^{y_1}_\lambda}
(f(u_\lambda)-f(u))(u_\lambda-u)^+dxdy.\label{eq yx2}
\end{eqnarray}
Now we only need to prove that
\begin{equation}\label{eq1}
 \int_{\Sigma^{y_1}_\lambda}(f(u_\lambda)-f(u))(u_\lambda-u)^+dxdy
<+\infty.
\end{equation}
In fact, for any given $\lambda\in\R$, using (\ref{eq1.2}), we choose  $R>1$  such that
$$0< u_\lambda(x,y)\leq C_\epsilon(1+|x|)^{-N-2\alpha} e^{-\theta_1|y_\lambda|} <s_0, \ \ \ \forall (x,y)\in B^c_R,$$
where $y_\lambda=(2\lambda-y_1,y')$ for $y=(y_1,y')\in\R^M$ and $s_0$ is from $(F)$.

If $u_\lambda(x,y)> u(x,y)$ for some $(x,y)\in{\Sigma^{y_1}_\lambda}\cap B^c_R$,  we have
$0<u(x,y)<u_\lambda(x,y)<s_0$. Using (\ref{eq rq2}) with $v=u_\lambda(x,y)$, then
$$ \frac{f(u_\lambda(x,y))-f(u(x,y))}{u_\lambda(x,y)-u(x,y)}
\leq \bar cu^\gamma_\lambda(x,y),$$
then
\begin{eqnarray*}
(f(u_\lambda(x,y))-f(u(x,y)))^+(u_\lambda(x,y)-u(x,y))^+
\leq \bar c u^{\gamma+2}_\lambda(x,y).
\end{eqnarray*}
The inequality above is obvious if $u_\lambda(x,y)\leq u(x,y)$ for some  $(x,y)\in{\Sigma^{y_1}_\lambda}\cap B^c_R$. Then
\begin{eqnarray*}
(f(u_\lambda)-f(u))^+(u_\lambda-u)^+\leq \bar c u^{\gamma+2}_\lambda \ \ \  {\rm in}\ \ {\Sigma^{y_1}_\lambda}\cap B^c_R.
 \end{eqnarray*}
Therefore,
\begin{eqnarray*}
&&\int_{{\Sigma^{y_1}_\lambda}\cap B^c_R}(f(u_\lambda)-f(u))^+(u_\lambda-u)^+dxdy
\\&\leq& \bar c\int_{{\Sigma^{y_1}_\lambda}\cap B^c_R}u^{\gamma+2}_\lambda(x,y)dxdy
\\&\leq&  \bar c C_\epsilon\int_{\Sigma^{y_1}_\lambda}(1+|x|)^{-(N+2\alpha)(\gamma+2)} e^{-(\gamma+2)\theta_1|y_\lambda|}dxdy
\\&\leq& \bar c C_\epsilon\int_{\R^N}(1+|x|)^{-(N+2\alpha)(\gamma+2)}dx\int_{\R^M} e^{-(\gamma+2)\theta_1|y|}dy
<+\infty,
\end{eqnarray*}
where the last inequality holds by $\gamma>\frac{2\alpha N}{(N+M)(N+2\alpha)}$. Since $u$ and $u_\lambda$ are bounded and $f$ is locally Lipschitz, we have
\begin{eqnarray*}
 \int_{\Sigma^{y_1}_\lambda\cap B_R}(f(u_\lambda)-f(u))^+(u_\lambda-u)^+dxdy
<+\infty.
\end{eqnarray*}
Therefore, (\ref{eq1}) holds. Together with (\ref{eq yx2}), we have the second inequality in the result.

{Next we show that the first inequality holds in
Lemma \ref{lemma dem4}.}
Let us denote
\begin{equation}\label{eq dem5} w(x,y)=\left\{
\arraycolsep=1pt
\begin{array}{lll}
(u_\lambda-u)^+(x,y),\ \ \ \ &
(x,y)\in \Sigma^{y_1}_\lambda,\\[2mm]
(u_\lambda-u)^-(x,y),\ \ \ \ & (x,y)\in (\Sigma^{y_1}_\lambda)^c
\end{array}
\right.
\end{equation}
and $${\rm{supp}}(w)=\overline{\{(x,y)\in\R^N\times\R^M\ |\
w(x,y)\not=0\}},$$
 where $(u_\lambda-u)^+(x,y)=\max\{(u_\lambda-u)(x,y),\ 0\}$,
$(u_\lambda-u)^-(x,y)=\min\{(u_\lambda-u)(x,y),\ 0\}$.
We observe that $w(x,y_1,y')=-w(x,2\lambda-y_1,y')$ for  $(x,y_1,y')\in\R^N\times\R\times\R^{M-1}$ and
\begin{equation}\label{eq 1.1}
w=u_\lambda-u\ \quad {\rm in}\ \ {\rm{supp}}(w).
\end{equation}
It is obvious that for $(x,y)\in\Sigma^{y_1}_\lambda\cap{\rm{supp}}(w)$,
$\{z\in\R^N| \  (z,y)\in(\Sigma^{y_1}_\lambda)^c\}=\O $ and
\begin{eqnarray*}
\R^N&=&\{z\in\R^N| \ (z,y)\in{\Sigma^{y_1}_\lambda\cap{\rm{supp}}(w)}\}\cup\\
&&\{z\in\R^N| \  (z,y)\in{\Sigma^{y_1}_\lambda\cap({\rm{supp}}(w))^c}\}\cup
\{z\in\R^N| \  (z,y)\in(\Sigma^{y_1}_\lambda)^c\}.
\end{eqnarray*}
Combining with (\ref{eq 1.1}), then for  $(x,y)\in\Sigma^{y_1}_\lambda\cap{\rm{supp}}(w)$,
\begin{eqnarray}
&&(-\Delta)_x^{\alpha}w(x,y)-(-\Delta)_x^{\alpha}(u_\lambda-u)(x,y)\nonumber
=\int_{\R^N}\frac{(u_\lambda-u)(z,y)-w(z,y)}{|x-z|^{N+2\alpha}}dz\nonumber
\\&&=\int_{\{z\in\R^N: (z,y)\in{\Sigma^{y_1}_\lambda\cap({\rm{supp}}(w))^c}\}}\frac{(u_\lambda-u)(z,y)}{|x-z|^{N+2\alpha}}dz
\leq 0,\label{desiwx2}
\end{eqnarray}
where the last inequality holds by  $u_\lambda-u\leq0$ in
$\Sigma^{y_1}_\lambda\cap({\rm{supp}}(w))^c$.
On one hand, from \equ{desiwx2} and $w=(u_\lambda -u)^+>0$ in $
\Sigma^{y_1}_\lambda\cap {\rm{supp}}(w)$, we have that
\begin{equation}\label{cotaD0}
 \int_{ \Sigma^{y_1}_\lambda\cap {\rm{supp}}(w)} (-\Delta)_x^{{\alpha}}w\,
 wdxdy
 \le  \int_{ \Sigma^{y_1}_\lambda\cap {\rm{supp}}(w)} (-\Delta)_x^{{\alpha}}(u_\lambda-u)(u_\lambda
 -u)^+dxdy.
  \end{equation}
On the other hand, we know that  $w(x,y)=(u_\lambda-u)(x,y)$ and
$(-\Delta)_yw(x,y)=(-\Delta)_y(u_\lambda-u)(x,y)$ for $(x,y)\in
\Sigma^{y_1}_\lambda\cap {\rm{supp}}(w)$. Together with
(\ref{cotaD0}), then
\begin{eqnarray*}
&&\int_{\Sigma^{y_1}_\lambda\cap{\rm{supp}}(w)}[(-\Delta)_x^{\alpha}w+(-\Delta)_yw+w]\,wdxdy
\\&\leq&\int_{\Sigma^{y_1}_\lambda\cap{\rm{supp}}(w)}[(-\Delta)_x^{\alpha}(u_\lambda-u)+(-\Delta)_y(u_\lambda-u)+(u_\lambda-u)](u_\lambda-u)^+dxdy
\end{eqnarray*}
and then by the fact of $w=(u_\lambda-u)^+=0$ in
$\Sigma^{y_1}_\lambda\cap({\rm{supp}}(w))^c$, we have that
\begin{eqnarray}
&&\int_{\Sigma^{y_1}_\lambda}[(-\Delta)_x^{\alpha}w+(-\Delta)_yw+w]\,wdxdy\nonumber
\\&\leq&\int_{\Sigma^{y_1}_\lambda}[(-\Delta)_x^{\alpha}(u_\lambda-u)+(-\Delta)_y(u_\lambda-u)+(u_\lambda-u)](u_\lambda-u)^+.\label{eq
2ap1}
\end{eqnarray}
By the definition of $w$, we have that
$$
\int_{\R^{N+M}}|w|^2dxdy=2\int_{\Sigma^{y_1}_\lambda}|w|^2dxdy,
$$
$$
\int_{\R^{N+M}}|w|^{\frac{2(N+M)}{N+M-2\alpha}}dxdy=2\int_{\Sigma^{y_1}_\lambda}|w|^{\frac{2(N+M)}{N+M-2\alpha}}dxdy,
$$
$$
\int_{\R^{N+M}}(-\Delta)_yw\,wdxdy=2\int_{\Sigma^{y_1}_\lambda}(-\Delta)_yw\,wdxdy,
$$
$$
\int_{\R^{N+M}}(-\Delta)_x^{\alpha}w\,wdxdy=2\int_{\Sigma^{y_1}_\lambda}(-\Delta)_x^{\alpha}w\,wdxdy,
$$
then, together with  Proposition \ref{lemma 1}, we obtain that
\begin{eqnarray}
&&\int_{\Sigma^{y_1}_\lambda}[(-\Delta)_x^{\alpha}w+(-\Delta)_yw+w]\,wdxdy\nonumber
\\&=&\frac{1}{2}\int_{\R^{N+M}}[(-\Delta)_x^{\alpha}w+(-\Delta)_yw+w]\,wdxdy\nonumber
\\&\geq& c_{3}(\int_{\R^{N+M}}|w|^{\frac{2(N+M)}{N+M-2\alpha}}dxdy)^{\frac{{N+M-2\alpha}}{N+M}}\nonumber
\\&=&c_{3}(2\int_{\Sigma^{y_1}_\lambda}|w|^{\frac{2(N+M)}{N+M-2\alpha}}dxdy)^{\frac{{N+M-2\alpha}}{N+M}},\label{eq
2ap}
\end{eqnarray}
for some $c_{3}>0$. Combining  (\ref{eq
2ap1}) with (\ref{eq 2ap}), by $w=(u_\lambda-u)^+$ in $\Sigma^{y_1}_\lambda$, we get the
first inequality  in Lemma \ref{lemma dem4}. The proof is
complete. \hfill$\Box$

\begin{lemma}\label{lemma dem411}
Under the assumptions of Theorem  \ref{teo r1jt}, for any
$\lambda\in\R$, there exists $c_{4}>0$ independent of
$\lambda$ such that
\begin{eqnarray*}
&&c_{4}(\int_{\Sigma^{x_1}_\lambda}|(u_\lambda-u)^+|^{\frac{2(N+M)}{N+M-2\alpha}}dxdy)^{\frac{N+M-2\alpha}{N+M}}
\\&\leq&
\int_{\Sigma^{x_1}_\lambda}[(-\Delta)_x^{\alpha}({u_\lambda}-u)+(-\Delta)_y({u_\lambda}-u)+({u_\lambda}-u)](u_\lambda-u)^+dxdy
<\infty,
\end{eqnarray*}
where $\Sigma^{x_1}_\lambda=\{(x_1,x',y)\in \R\times\R^{N-1}\times\R^M\  |\  x_1>\lambda\}.$
\end{lemma}
{\bf Proof.} The proof proceeds similarly to the proof of Lemma \ref{lemma dem4}, the only difference is
to show (\ref{desiwx2}) with $(x,y)\in
\Sigma^{x_1}_\lambda\cap{\rm{supp}}(w)$. It is obvious that
\begin{eqnarray*}
\R^N&=&\{z\in\R^N| \ (z,y)\in{\Sigma^{x_1}_\lambda\cap{\rm{supp}}(w)}\}\cup\\
&&\{z\in\R^N| \  (z,y)\in{\Sigma^{x_1}_\lambda\cap({\rm{supp}}(w))^c}\}\cup\\
&&\{z\in\R^N| \  (z,y)\in(\Sigma^{x_1}_\lambda)^c\cap({\rm{supp}}(w))^c\}\cup\\
&&\{z\in\R^N| \  (z,y)\in(\Sigma^{x_1}_\lambda)^c\cap{\rm{supp}}(w)\}
\end{eqnarray*}
and $w=u_\lambda-u$ in ${\rm{supp}}(w)$, then for  $(x,y)\in\Sigma^{x_1}_\lambda\cap{\rm{supp}}(w)$,
\begin{eqnarray*}
&&(-\Delta)_x^{\alpha}w(x,y)-(-\Delta)_x^{\alpha}(u_\lambda-u)(x,y)\nonumber
=\int_{\R^N}\frac{(u_\lambda-u)(z,y)-w(z,y)}{|x-z|^{N+2\alpha}}dz\nonumber
\\&=&\int_{\{z\in\R^N | \ (z,y)\in{\Sigma^{x_1}_\lambda\cap({\rm{supp}}(w))^c}\}}(\frac{1}{|x-z|^{N+2\alpha}}-\frac{1}{|x-z_\lambda|^{N+2\alpha}})(u_\lambda-u)(z,y)dz
\\&\leq& 0,
\end{eqnarray*}
where $z_\lambda=(2\lambda-z_1,z')$ for $z=(z_1,z')\in \R^N$ and
 the last inequality holds by   $u_\lambda-u\leq0$ in
$\Sigma^{x_1}_\lambda\cap({\rm{supp}}(w))^c$.
\hfill$\Box$

\begin{teo}\label{teo w1}
Under the assumptions of Theorem  \ref{teo r1jt}, for $x\in\R^N$, we have
$$u(x,y)=u(x,|y|)$$
and $u$ is strictly decreasing in $y$-direction.
\end{teo}
{\bf Proof.}
We divide the proof into three steps.

\noindent{\bf Step 1:} $\lambda_0:=\sup\{\lambda\ |\ u_\lambda\leq u\  {\rm in}\ \Sigma^{y_1}_\lambda \}$ is finite.
Since $u$ decays at infinity, we observe that the set $\{\lambda\ |\ u_\lambda\leq u\  {\rm in}\ \Sigma^{y_1}_\lambda \}$ is nonempty.
Using $(u_\lambda-u)^+$ as a test function in the equation for
$u$ and $u_\lambda$,  by
(\ref{eq rq2}) and  H\"{o}lder inequality, for $\lambda$ big (negative), we find that
\begin{eqnarray*}
&&\int_{\Sigma^{y_1}_\lambda}[(-\Delta)_x^{\alpha}({u_\lambda}-u)+(-\Delta)_y({u_\lambda}-u)+({u_\lambda}-u)](u_\lambda-u)^+dxdy
\\&&= \int_{\Sigma^{y_1}_\lambda} (f(u_\lambda)-f(u))(u_\lambda-u)^+dxdy
\\&&= \int_{\Sigma^{y_1}_\lambda}\frac{f(u_\lambda)-f(u)}{u_\lambda-u}[(u_\lambda-u)^+]^2dxdy
\leq \bar c\int_{\Sigma^{y_1}_\lambda}{u^\gamma_\lambda}[(u_\lambda-u)^+]^2dxdy
\\&&\leq c_{5}\int_{\Sigma^{y_1}_\lambda}(1+|x|)^{-\gamma(N+2\alpha)} e^{-\gamma\theta_1|y_\lambda|}[(u_\lambda-u)^+]^2dxdy
\\&&\leq  c_{5}(\int_{\Sigma^{y_1}_\lambda}(1+|x|)^{-a} e^{-b|y_\lambda|}dxdy)^{\frac{2\alpha}{N+M}}
  (\int_{\Sigma^{y_1}_\lambda}|(u_\lambda-u)^+|^{\frac{2(N+M)}{N+M-2\alpha}}dxdy)^{\frac{N+M-2\alpha}{N+M}},
\end{eqnarray*}
where $a=\frac{\gamma(N+2\alpha)(N+M)}{2\alpha}$ and $b=\frac{\theta_1\gamma(N+M)}{2\alpha}$.
Since $\gamma>\frac{2\alpha N}{(N+2\alpha)(N+M)}$, we have that $a>N$. Then we can choose $R>0$ such that for all $\lambda<-R$,
$$c_{5}(\int_{\Sigma_\lambda^{y_1}}(1+|x|)^{-a} e^{-b|y_\lambda|}dxdy)^{\frac{2\alpha}{N+M}}\leq\frac{1}{4}.$$
By  Lemma \ref{lemma dem4}, we obtain that
\begin{eqnarray*}
\int_{\Sigma^{y_1}_\lambda}|(u_\lambda-u)^+|^{\frac{2(N+M)}{N+M-2\alpha}}dxdy=0,\ \ \quad \forall\ \lambda<-R.
\end{eqnarray*}
Thus $u_\lambda\leq u$ in $\Sigma^{y_1}_\lambda$ for all $\lambda<-R$ and then   conclude that $\lambda_0\geq -R.$
On the other hand, since $u$ decays at infinity, then there exist $\lambda_1\in\R$ and $(x,y)\in \Sigma^{y_1}_\lambda$ such that
$u(x,y)<u_{\lambda_1}(x,y)$.
 Hence  $\lambda_0$ is finite.
 \medskip

\noindent{\bf Step 2: }$u\equiv u_{\lambda_0}$ in $\Sigma_{\lambda_0}^{y_1}$.
Assuming the contrary, we have that $u\not\equiv u_{\lambda_0}$ and  $u\geq
u_{\lambda_0}$ in $\Sigma_{\lambda_0}^{y_1}$, in this case the following
claim holds.
\smallskip

 \noindent{\bf Claim 1.} \emph{If $u\not\equiv
u_{\lambda_0}$ and $u\geq u_{\lambda_0}$ in $\Sigma_{\lambda_0}^{y_1}$,
then  $u>u_{\lambda_0}$ in $\Sigma_{\lambda_0}^{y_1}$.}

Let us assume, for the moment, that Claim 1 is true, then for any given
$\lambda\in({\lambda_0},{\lambda_0}+\epsilon)$, where $\epsilon>0$ is chosen later.
Let $P=(0,\cdots,\lambda,\cdots,0)\in T_\lambda^{y_1}$ and $B(P,R)$ be the ball centered at $P$ and with
radius $R>1$ to be chosen later.
 Define $B_1=\Sigma_\lambda^{y_1}\cap B(P,R)$ and let us consider  $(u_\lambda-u)^+$  test function in the equation for
$u$ and $u_\lambda$ in $\Sigma_\lambda^{y_1}$, then from Lemma \ref{lemma dem4}
we obtain
\begin{eqnarray}
&&(\int_{\Sigma^{y_1}_\lambda}|(u_\lambda-u)^+|^{\frac{2(N+M)}{N+M-2\alpha}}dxdy)^{\frac{N+M-2\alpha}{N+M}}\nonumber
\\&\leq& c_{6}\int_{\Sigma^{y_1}_\lambda}[(-\Delta)_x^{\alpha}({u_\lambda}-u)+(-\Delta)_y({u_\lambda}-u)+({u_\lambda}-u)](u_\lambda-u)^+dxdy\nonumber
\\&=&  c_{6}\int_{\Sigma_\lambda ^{y_1}} (f(u_\lambda)-f(u))(u_\lambda-u)^+dxdy.\label{E0}
\end{eqnarray}
We estimate the integral on the right. Proceeding  as in Step 1, we can choose $R>1$ big enough such that
$$ c_{7}(\int_{\Sigma_\lambda^{y_1}\setminus{B_1}}(1+|x|)^{-a} e^{-b|y_\lambda|}dxdy)^{\frac{2\alpha}{N+M}}\leq\frac{1}{4}$$
for some $c_{7}>0$,
where $a=\frac{\gamma(N+2\alpha)(N+M)}{2\alpha}$ and $b=\frac{\theta_1\gamma(N+M)}{2\alpha}$.
Then
\begin{eqnarray}
&&\int_{\Sigma_\lambda^{y_1}\setminus{B_1}}(f(u_\lambda)-f(u))(u_\lambda-u)^+dxdy
\leq \bar c \int_{\Sigma_\lambda^{y_1}\setminus{B_1}}u^\gamma_\lambda |(u_\lambda-u)^+|^2dxdy\nonumber
\\&&\leq c_{7}(\int_{\Sigma_\lambda^{y_1}\setminus{B_1}}(1+|x|)^{-a} e^{-b|y_\lambda|}dxdy)^{\frac{2\alpha}{N+M}}
  (\int_{\Sigma^{y_1}_\lambda}|(u_\lambda-u)^+|^{\frac{2(N+M)}{N+M-2\alpha}}dxdy)^{\frac{N+M-2\alpha}{N+M}}\nonumber
\\&&\le\frac{1}{4} (\int_{\Sigma_\lambda^{y_1}}|(u_\lambda-u)^+|^{\frac{2(N+M)}{N+M-2\alpha}}dxdy)^{\frac{N+M-2\alpha}{N+M}}.
\label{E2}\end{eqnarray}
Now using Claim 1,  we  choose $\epsilon>0$ such that
$c_{8}|B_1\cap{{\rm{supp}}{(u_\lambda-u)^+}}|^{\frac{2\alpha}{N+M}}<1/4,$
for some $c_{8}>0$.
Since $f$ is locally Lipschitz, using H\"older inequality, we have
\begin{eqnarray}
 && \int_{B_1} (f(u_\lambda)-f(u))(u_\lambda-u)^+dxdy\leq c_{46}\int_{B_1}|(u_\lambda-u)^+|^2\chi_{{\rm{supp}}{(u_\lambda-u)^+}}dxdy\nonumber
\\&&= c_{8}|B_1\cap{{\rm{supp}}{(u_\lambda-u)^+}}|^{\frac{2\alpha}{N+M}}
  (\int_{B_1}|(u_\lambda-u)^+|^{\frac{2(N+M)}{N+M-2\alpha}}dxdy)^{\frac{N+M-2\alpha}{N+M}}\nonumber
\\&&\leq\frac{1}{4} (\int_{B_1}|(u_\lambda-u)^+|^{\frac{2(N+M)}{N+M-2\alpha}}dxdy)^{\frac{N+M-2\alpha}{N+M}}.\label{E1}
\end{eqnarray}
From \equ{E0}, \equ{E2} and \equ{E1},
 it follows that
 $(u_\lambda-u)^+=0$ in $\Sigma_\lambda^{y_1}$. Then $ u_{\lambda}\leq u$ in
$\Sigma_\lambda^{y_1}$ for $\lambda\in({\lambda_0},{\lambda_0}+\epsilon)$, which contradicts the definition of $\lambda_0$.
As a consequence, we have $u\equiv u_{\lambda_0}$ in $\Sigma_{\lambda_0}^{y_1}$.

In order to complete Step 2, we only need to prove Claim 1.
\smallskip

\noindent{\bf Proof of Claim 1.}
By contradiction, if  there exists $(\bar x,\bar y)\in \Sigma_{\lambda_0}^{y_1}$ such that
$u(\bar x,\bar y)=u_{\lambda_0}(\bar x,\bar y),$ then
\begin{eqnarray*}
&&(-\Delta)_x^{\alpha} (u-u_{\lambda_0})(\bar x,\bar y)+(-\Delta)_y
(u-u_{\lambda_0})(\bar x,\bar y)+(u-u_{\lambda_0})(\bar x,\bar y)
\\&&=f(u(\bar x,\bar y))-f(u_{\lambda_0}(\bar x,\bar y))=0.
\end{eqnarray*}
Since $(u-u_{\lambda_0})(\bar x,\bar y)=\min_{\Sigma_{\lambda_0}^{y_1}}(u-u_{\lambda_0})=0$, we have
$(-\Delta)_y(u-u_{\lambda_0})(\bar x,\bar y)\leq 0$, then
\begin{equation}\label{eq wy1}
(-\Delta)_x^{\alpha} (u-u_{\lambda_0})(\bar x,\bar y)\geq 0.
\end{equation}
The other side, we observe that  $\{z\in\R^N | \   (z,\bar y)\in(\Sigma_{\lambda_0}^{y_1})^c\}=\O$ when $(\bar x,\bar y)\in \Sigma_{\lambda_0}^{y_1}$.
By $u(\bar x,\bar y)=u_{\lambda_0}(\bar x,\bar y)$ and then
\begin{eqnarray}
&&(-\Delta)_x^{\alpha} (u-u_{\lambda_0})(\bar x,\bar y)
=-\int_{\R^N}\frac{(u-u_{\lambda_0})(z,\bar y)}{|\bar x-z|^{N+2\alpha}}dz\nonumber
\\&&=-\int_{\{z\in\R^N| \ (z,\bar y)\in\Sigma_{\lambda_0}^{y_1}\}}\frac{(u-u_{\lambda_0})(z,\bar y)}{|\bar x-z|^{N+2\alpha}}dz
\leq0,\label{eq wy2}
\end{eqnarray}
where the last inequality holds by
$u\geq u_{\lambda_0}$ in $\Sigma_{\lambda_0}^{y_1}$.

Combining (\ref{eq wy1}) with (\ref{eq wy2}), we obtain that $(-\Delta)_x^{\alpha} (u-u_{\lambda_0})(\bar x,\bar y)=0$ and then from (\ref{eq wy2}),
we have that
\begin{equation}\label{eq 1.20}
u(z,\bar y)=u_{\lambda_0}(z,\bar y),\quad \forall z\in\R^N,
\end{equation}
 this means that $u-u_{\lambda_0}$ has property $(P)$ and  by $u\neq u_{\lambda_0}$ in $\Sigma_{\lambda_0}^{y_1}$ we  have
$$(\bar x, \bar y)\in(\Sigma_{\lambda_0}^{y_1})_0:=\{(x,y)\in\Sigma_{\lambda_0}^{y_1}\ | \ (u-u_{\lambda_0})(x,y)=\inf_{\Sigma_{\lambda_0}^{y_1}}(u-u_{\lambda_0})=0\}\subsetneqq \Sigma_{\lambda_0}^{y_1}.$$
Moreover, by Proposition \ref{proposition 2.1} with $\Omega=\Sigma_{\lambda_0}^{y_1}$, we observe that $\Sigma_{\lambda_0}^{y_1}\setminus (\Sigma_{\lambda_0}^{y_1})_0$
satisfies interior cylinder condition at point $(x_0,y_0)\in \partial(\Sigma_{\lambda_0}^{y_1})_0\cap \Sigma_{\lambda_0}^{y_1}$.
Then there exist $r>0$ small and $\tilde y\in \R^M$ such that
$$ O_r:=B_r^N(x_0)\times B_r^M(\tilde y)\subset\Sigma_{\lambda_0}^{y_1}\setminus (\Sigma_{\lambda_0}^{y_1})_0\quad
{\rm{and}}\quad (x_0,y_0)\in\partial O_r.$$
Let  $D$ be defined by (\ref{eq ms1}). Since $u\geq u_{\lambda_0}$ in $\Sigma_{\lambda_0}^{y_1}$,
then for any $(x,y)\in D$, we have
\begin{eqnarray*}
\int_{\R^N\setminus B_r^N(x_0)}\frac{(u-u_{\lambda_0})(z,y)}{|x-z|^{N+2\alpha}}dz\ge0.
\end{eqnarray*}
Finally, it is obvious that
 $$(-\Delta)_x^{\alpha} (u-u_{\lambda_0})+(-\Delta)_y (u-u_{\lambda_0})+h(u-u_{\lambda_0})=0\quad {\rm{in}} \ \Sigma_{\lambda_0}^{y_1},$$
where $h=1-\frac{f(u)-f(u_{\lambda_0})}{u-u_{\lambda_0}}\in L_{loc}^\infty(\Sigma_{\lambda_0}^{y_1})$.
Then we use Theorem \ref{teo SMP} to obtain
$$u\equiv u_{\lambda_0} \quad {\rm{in}}\ \tilde\Sigma_{\lambda_0}^{y_1},$$
which contradicts the condition of $u\neq u_{\lambda_0}$  in $\Sigma_{\lambda_0}^{y_1}$,
then we obtain the results in Claim 1.
 \medskip

\noindent\textbf{Step 3.} By translation, we may say that
$\lambda_0=0.$ Repeating the argument from the other side, we
find that $u$ is symmetric about $y_1$-axis. Using the same argument
in any $y$-direction, we  conclude that
$$u(x,y)=u(x,|y|),\quad (x,y)\in\R^N\times\R^M.$$

Finally, {we prove that $u(x,|y|)$ is strictly decreasing in  $|y|>0$.}
Indeed, for any given  $y_1<\widetilde{y}_1<0$
and letting  $\lambda=\frac{y_1+\widetilde{y}_1}{2}$. Then,  as proved above  we have
$$u>u_\lambda\quad \  {\rm{in}} \quad \Sigma_{\lambda}^{y_1}.$$
For any given $x\in \R^N$, we observe that $(x,\widetilde y_1,0,\cdots,0)\in \Sigma_{\lambda}^{y_1}$, then
\begin{eqnarray*}
u(x,\widetilde y_1,0,\cdots,0)>u_\lambda(x,\widetilde y_1,0,\cdots,0)=u(x,y_1,0,\cdots,0).
\end{eqnarray*}
 Using the result of $u(x,y)=u(x,|y|)$ for all $(x,y)\in\R^N\times\R^M$ and $|\widetilde{y}_1|<|y_1|$, we conclude  monotonicity of $u$ respect to $y$.
This completes the proof.
 \hfill$\Box$

\medskip

Next we study the symmetry result in $x$-direction.

\begin{teo}\label{teo w2}
Under the assumptions of Theorem  \ref{teo r1jt}, for $y\in\R^M$, we have
$$u(x,y)=u(|x|,y)$$
and $u$ is strictly decreasing in $x$-direction.
\end{teo}
{\bf Proof.}  The proof of this theorem goes like the one for Theorem \ref{teo w1}.
The only place where there is a difference is in the following property:
\emph{if $u\not\equiv u_{\lambda_0}$ and $u\geq u_{\lambda_0}$ in $\Sigma_{\lambda_0}^{x_1}$,
then  $u>u_{\lambda_0}$ in $\Sigma_{\lambda_0}^{x_1}$.}
By contradiction, if there exists  $(\bar x,\bar y)\in \Sigma_{\lambda_0}^{x_1}$ such that $u(\bar x,\bar y)=u_{\lambda_0}(\bar x,\bar y)$, then
\begin{eqnarray*}
&&(-\Delta)_x^{\alpha} (u-u_{\lambda_0})(\bar x,\bar y)+(-\Delta)_y
(u-u_{\lambda_0})(\bar x,\bar y)+(u-u_{\lambda_0})(\bar x,\bar y)
\\&&=f(u(\bar x,\bar y))-f(u_{\lambda_0}(\bar x,\bar y))=0.
\end{eqnarray*}
Since $u\geq u_{\lambda_0}$ in $\Sigma_{\lambda_0}^{x_1}$, we have  $(u-u_{\lambda_0})(\bar x,\bar y)=\min_{\Sigma_{\lambda_0}^{x_1}}(u-u_{\lambda_0})=0$ and
$(-\Delta)_y(u-u_{\lambda_0})(\bar x,\bar y)\leq 0$ and then
$$
(-\Delta)_x^{\alpha} (u-u_{\lambda_0})(\bar x,\bar y)\geq 0.
$$
The other side, by direct computation, we have that
\begin{eqnarray*}
&&(-\Delta)_x^{\alpha} (u-u_{\lambda_0})(\bar x,\bar y)
=\int_{\R^N}\frac{(u_{\lambda_0}-u)(z,\bar y)}{|\bar x-z|^{N+2\alpha}}dz\nonumber
\\&&=\int_{\{z\in\R^N| \ (z,\bar y)\in\Sigma_{\lambda_0}^{x_1}\}}(\frac{1}{|\bar x-z|^{N+2\alpha}}-\frac{1}{|\bar x-z_{\lambda_0}|^{N+2\alpha}})
(u_{\lambda_0}-u)(z,\bar y)dz
\leq0,
\end{eqnarray*}
where $z_{\lambda_0}=(2{\lambda_0}-z_1,z')$ for $z=(z_1,z')\in \R^N$ and  the last inequality holds by
$u\geq u_{\lambda_0}$ in $\Sigma_{\lambda_0}^{x_1}$. Therefore,
\begin{equation}\label{eq 1.206}
u(z,\bar y)=u_{\lambda_0}(z,\bar y),\quad \forall z\in\R^N,
\end{equation}
 this means that $u-u_{\lambda_0}$ has property $(P)$ and  by $u\neq u_{\lambda_0}$ in $\Sigma_{\lambda_0}^{x_1}$ we  have
$$(\bar x, \bar y)\in(\Sigma_{\lambda_0}^{x_1})_0:=\{(x,y)\in\Sigma_{\lambda_0}^{x_1}\ | \ (u-u_{\lambda_0})(x,y)=\inf_{\Sigma_{\lambda_0}^{x_1}}(u-u_{\lambda_0})=0\}\subsetneqq \Sigma_{\lambda_0}^{x_1}.$$
Moreover, by Proposition \ref{proposition 2.1}, we observe that $\Sigma_{\lambda_0}^{x_1}\setminus (\Sigma_{\lambda_0}^{x_1})_0$
satisfies interior cylinder condition at point $(x_0,y_0)\in \partial(\Sigma_{\lambda_0}^{x_1})_0\cap \Sigma_{\lambda_0}^{x_1}$.
Then there exist $r_1>0$ and $\tilde y\in \R^M$ such that for all $r\in(0,r_1]$,
$$ O_r:=B_r^N(x_0)\times B_r^M(\tilde y)\subset\Sigma_{\lambda_0}^{x_1}\setminus (\Sigma_{\lambda_0}^{x_1})_0\quad
{\rm{and}}\quad (x_0,y_0)\in\partial O_r.$$
 Next we show that there exists some $r\in(0,r_1]$ such that for any $(x,y)\in D$,
\begin{equation}\label{eq 2.8}
\int_{\R^N\setminus B_r^N(x_0)}\frac{(u-u_{\lambda_0})(z,y)}{|x-z|^{N+2\alpha}}dz\ge0,
\end{equation}
where  $D$ is defined by (\ref{eq ms1}).
Indeed, since $u\not\equiv u_{\lambda_0}$ and $u\geq u_{\lambda_0}$ in $\Sigma_{\lambda_0}^{x_1}$,
then for $(x,y)\in D\subset\Sigma_{\lambda_0}^{x_1}$,
we have that
$$\int_{\R^N}\frac{(u-u_{\lambda_0})(z,y)}{|x-z|^{N+2\alpha}}dz>0.$$
Let us define
\begin{equation}\label{eq 2.82}
r(x,y)=\sup\{r\in(0,r_1]: \int_{\R^N\setminus B_r^N(x_0)}\frac{(u-u_{\lambda_0})(z,y)}{|x-z|^{N+2\alpha}}dz\ge0\}.
\end{equation}
Let $r_m=\inf_{(x,y)\in D}r(x,y)$, it is obvious that $r_m\in[0,r_1]$. Now we prove that $r_m>0$. By contradiction, if $r_m=0$,
then there exist a sequence $(x_n,y_n)\in D$ and $(\tilde{x},\tilde{y})\in \bar D$ such that
$(x_n,y_n)\to (\tilde{x},\tilde{y})$ and $r(x_n,y_n)\to 0,$ as $n\to +\infty.$
Since $r(x,y)$ is continuous, then $r(\tilde{x},\tilde{y})=0.$
If $(\tilde{x},\tilde{y})\in\bar D\setminus (\Sigma_{\lambda_0}^{x_1})_0$, i.e. $u(\tilde{x},\tilde{y})>u_{\lambda_0}(\tilde{x},\tilde{y})$,
we have
\begin{eqnarray*}
&&\int_{\R^N}\frac{(u-u_{\lambda_0})(z,\tilde{y})}{|\tilde{x}-z|^{N+2\alpha}}dz
\\&=&\int_{\{z\in\R^N| \ (z, \tilde{y})\in\Sigma_{\lambda_0}^{x_1}\}}(u-u_{\lambda_0})(z,\tilde{y})(\frac{1}{|\tilde{x}-z|^{N+2\alpha}}
-\frac{1}{| \tilde{x}-z_{\lambda_0}|^{N+2\alpha}})
dz
>0.
\end{eqnarray*}
By the continuity of the integration and (\ref{eq 2.82}), we obtain that $r(\tilde{x},\tilde{y})>0$, which is impossible.

Then $(\tilde{x},\tilde{y})\in\bar D\cap(\Sigma_{\lambda_0}^{x_1})_0$, i.e. $u(\tilde{x},\tilde{y})=u_{\lambda_0}(\tilde{x},\tilde{y})$.
Since the function $u-u_{\lambda_0}$ has property $(P)$, then for any $\tilde{r}>0$,
$$\int_{\R^N\setminus B_{\tilde{r}}^N(x_0)}\frac{(u-u_{\lambda_0})(z,\tilde{y})}{|\tilde{x}-z|^{N+2\alpha}}dz=0.$$
Combining with (\ref{eq 2.82}), we obtain that $r(\tilde{x},\tilde{y})=r_1>0$, which contradicts $r(\tilde{x},\tilde{y})=0$.
As a consequence, we have that $0<r_m\leq r_1$. Taking $r=r_m$, then (\ref{eq 2.8}) holds for any $(x,y)\in D$.
Finally, it is obvious that
 $$(-\Delta)_x^{\alpha} (u-u_{\lambda_0})+(-\Delta)_y (u-u_{\lambda_0})+h(u-u_{\lambda_0})=0\quad {\rm{in}} \ \Sigma_{\lambda_0}^{x_1},$$
where $h=1-\frac{f(u)-f(u_{\lambda_0})}{u-u_{\lambda_0}}\in L_{loc}^\infty(\Sigma_{\lambda_0}^{x_1})$.
Then we use Theorem \ref{teo SMP} to obtain that
$$u\equiv u_{\lambda_0} \quad {\rm{in}}\ \tilde\Sigma_{\lambda_0}^{x_1},$$
which contradicts the condition of $u\neq u_{\lambda_0}$  in $\Sigma_{\lambda_0}^{x_1}$. Then  $u> u_{\lambda_0}$  in $\Sigma_{\lambda_0}^{x_1}$,
to complete the proof.
 \hfill$\Box$\\

\noindent {\bf Acknowledgements:}  P.F. was  partially supported by Fondecyt Grant \# 1110291 and
 BASAL-CMM projects.
Y.W. was partially supported by Becas CMM.

\end{document}